\documentclass[article,3p,12pt]{elsarticle}

\usepackage{color}
\usepackage{mathtools}
\usepackage{hyperref} 
\usepackage{setspace}
\usepackage{amsfonts}
\usepackage{amssymb}
\usepackage{amsthm}
\usepackage{bm}
\usepackage{ulem}
\usepackage{caption}
\usepackage{subfigure}
\usepackage{multirow}
\usepackage{array}
\usepackage{float}
\usepackage{todonotes}
\usepackage{textcomp}
\usepackage{pgfplots}
\pgfplotsset{width=10cm,compat=1.9}
\usepackage{adjustbox}
\usepackage{booktabs}

\usepackage{tikz}
\usetikzlibrary{calc}
\usetikzlibrary{arrows.meta, positioning, shapes.geometric, fit}
\tikzset{arrow/.style={-{Latex[length=3mm, width=1.5mm]}}}

\usepackage{txfonts}

\def\nobar{}

\numberwithin{equation}{section}

\makeatletter
\renewcommand{\@biblabel}[1]{#1\hfill \hspace{-0.2cm}}
\makeatother

\journal{ }

\begin{document}
\begin{frontmatter}
\title{\textsc{Mathematical modeling and sensitivity analysis of hypoxia-activated drugs}}

\author[dei]{Alessandro Coclite}

\author[mox]{Riccardo Montanelli Eccher}

\author[it]{Luca Possenti}

\author[mox]{Piermario Vitullo}

\author[mox]{Paolo Zunino\corref{cor}}
\ead{paolo.zunino@polimi.it}

\cortext[cor]{Corresponding author}

\address[dei]{Dipartimento di Ingegneria Elettrica e dell'Informazione (DEI), Politecnico di Bari, Bari, Italy}

\address[mox]{MOX, Dipartimento di Matematica, Politecnico di Milano, Milano, Italy}

\address[it]{Data Science Unit, Fondazione IRCCS Istituto Nazionale dei Tumori, Milano, Italy}

\begin{abstract}

Hypoxia-activated prodrugs offer a promising strategy for targeting oxygen-deficient regions in solid tumors, which are often resistant to conventional therapies. However, modeling their behavior is challenging because of the complex interplay between oxygen availability, drug activation, and cell survival. In this work, we develop a multiscale and mixed-dimensional model that couples spatially resolved drug and oxygen transport with pharmacokinetics and pharmacodynamics to simulate the cellular response. The model integrates blood flow, oxygen diffusion, and consumption, drug delivery, and metabolism. To reduce computational cost, we mitigate the global nonlinearity by coupling the multiscale and mixed-dimensional models one-way with a reduced 0D model for drug metabolization. The global sensitivity analysis is then used to identify key parameters influencing drug activation and therapeutic outcome. This approach enables efficient simulation and supports the design of optimized therapies targeting hypoxia.

\begin{keyword}
Mathematical oncology; Drug transport and metabolism; Hypoxia-activated prodrugs; Multiscale modeling; Numerical simulation; Global sensitivity analysis.

\textbf{Mathematics Subject Classification:} 92C50, 92C45, 35Q92, 35R20, 65M60.
\end{keyword}
\end{abstract}

\end{frontmatter}

\section*{Introduction}

Delivering effective chemo- or radiotherapy to solid tumors remains a significant challenge in oncology, as it is intricately linked to the tumor microenvironment (TME) and the physicochemical properties governing drug transport~\cite{Barker2015,dewhirst2017transport}. 
As highlighted in a comprehensive review~\cite{dewhirst2017transport}, the heterogeneous nature of tumor vasculature and tissue characteristics results in highly variable drug and oxygen distribution, necessitating a deeper understanding of the underlying transport processes and cellular responses. 
Mathematical modeling, coupled with spatially resolved experimental observations, is increasingly recognized as a crucial tool for advancing this understanding and designing improved treatment strategies. In this broader context, fractional-order approaches have also been successfully applied in mathematical epidemiology to capture memory and hereditary properties in biological systems~\cite{abdo2020, kumar2020, khan2022, pandey2021, alderremy2021}. Fractional-order models have recently gained attention in biomedical applications, as they can effectively capture memory and hereditary effects in complex biological systems. For instance, recent studies have investigated the qualitative analysis of generalized mixed fractional differential equations with applications to medicine~\cite{Hattaf2025}, and the dynamics of a reaction--diffusion fractional-order model for oncolytic virotherapy with CTL immune response~\cite{Chaos2022}. Recent mathematical models have also investigated the dynamics of cancer–immune interactions under therapeutic interventions. For instance, Ghosh et al.~\cite{Ghosh2025} analyzed HIV replication kinetics under dual inhibitors, providing insights into treatment efficacy through fractional modeling; Ali et al.~\cite{Ali2025} developed a mathematical framework for chemo-iPSC therapy in cervical cancer; and Kumar et al.~\cite{Kumar2019} studied combination drug therapy strategies for HPV-induced cervical cancer. These studies highlight the growing relevance of fractional-order and optimal-control approaches in describing complex tumor–immune–therapy systems.

A key factor influencing the success of cancer therapies, including chemotherapy and radiotherapy, is the presence of \textit{hypoxia}, low oxygen levels, within solid tumors~\cite{Barker2015, dewhirst2017transport, Hicks2006}. 
Hypoxia is a well-documented characteristic of many aggressive and immunosuppressive tumors, associated with genomic instability, apoptosis, angiogenesis, metastasis, invasion, and metabolic reprogramming~\cite{Li2018, Fu2021}. 
Furthermore, tumor hypoxia is a critical contributor to resistance to both radiotherapy and chemotherapy~\cite{Barker2015, Forster2019, Possenti2024}. 
In particular, oxygen is one of the most potent sensitizers in radiotherapy, since oxygen enhances the cytotoxic effect of radiation by stabilizing DNA-damaging free radicals~\cite{Barker2015}.
To specifically target these oxygen-deficient areas, \textit{hypoxia-activated prodrugs} (HAPs) have emerged as a promising therapeutic strategy~\cite{Sausville2006, Pruijn2008, Li2021}. 
These agents are designed to be preferentially activated under low oxygen conditions, ideally releasing cytotoxic compounds within hypoxic tumor regions while minimizing systemic toxicity.
However, modeling the behavior of hypoxia-activated drugs presents unique challenges due to the \textit{two-way coupling} between drug activation and oxygen availability, mediated by the survival fraction of cancer cells. 
Drug activation depends on the hypoxic state, which, in turn, is influenced by oxygen delivery through the microvasculature and oxygen consumption by viable tumor cells. 
As the drug exerts its cytotoxic effect, the number of viable cells decreases, potentially altering oxygen consumption rates, thus affecting the conditions that govern further drug activation~\cite{Linninger2024}. 
This intricate feedback loop requires sophisticated modeling approaches to capture these dynamic interactions. 

Computational models are essential for understanding the complex dynamics of the vascular and tumor microenvironment, where mechanical, hemodynamic, and regulatory factors interact in both healthy and diseased tissues. 
These models allow for the quantitative analysis of key processes such as blood flow, oxygen transport, and drug delivery, which are particularly relevant in cancer and radiotherapy, where they help assess oxygen dependence and vascular function to optimize treatment strategies.
Possenti and collaborators have developed advanced multiscale models to study tumor oxygenation and the impact of the vascular network on radiotherapy. Their work includes modeling microcirculation, fluid exchange, oxygen and drug transport, and applying global sensitivity analysis to identify key parameters. These tools support the exploration of hypothetical and \textit{in vitro} scenarios, improve our understanding of tumor responses, and help plan treatment~\cite{Cattaneo2014b, CZ, Possenti2019, Possenti2019b, Possenti2021}.
However, modeling complex scenarios such as hypoxia-activated drugs remains computationally demanding, as it involves coupling blood flow, interstitial dynamics, oxygen metabolism, and pharmacokinetics/pharmacodynamics across multiple spatial and temporal scales. 
Capturing these interactions with high spatial resolution can be extremely resource-intensive.
In the case of hypoxia-activated drugs, where drug activation depends on intricate interactions between blood flow, oxygen transport, and metabolic processes, understanding which parameters most influence treatment outcomes, such as drug concentration, oxygen levels, or cell survival, is crucial. However, techniques like sensitivity analysis, while informative, add computational cost. 
Sensitivity analysis is essential for extracting meaningful information from complex drug delivery models in the tumor microenvironment~\cite{Possenti2020, Vitullo2023}, but it requires significant computational demands. 
Although methods like Sobol's indices provide rigorous and quantitative assessments of the influence of the parameters, they require numerous simulations and are computationally intensive~\cite{saltelli2010variance}. 
More efficient approaches, such as the Morris elementary effects method, offer qualitative insight with fewer simulations, but still require careful sampling of the input space~\cite{Saltelli2008}. 
These challenges are amplified in models that include multiscale or multiphysics features, where each simulation can already be computationally costly. 
Thus, while powerful, these models often require simplification, such as the reduced-order model addressed in~\cite{Vitullo2024} or high-performance computing for practical use.

To address this complexity from a mathematical and computational point of view, we propose a \textit{multiscale} modeling approach that combines spatially distributed models of drug delivery in three dimensions (3D) with surrogate models (where spatial dependence is neglected). 
Spatially distributed models are essential to capture the heterogeneous oxygen distribution within the TME, influenced by the morphology of the microvascular network and the dynamics of oxygen transport. 
These models can account for blood flow in the vasculature (often represented as a 1D network embedded in a 3D tissue domain) and the transport and diffusion of oxygen and drugs within porous tumor tissue. 
Several research studies have focused on developing such a microcirculation model to study oxygen delivery and drug transport in tumors, highlighting the importance of vascular architecture and transport properties. 
In contrast, while neglecting spatial variations, surrogate models can provide computational efficiency in simulating cellular-level drug activation and its impact on cell survival over time.
The integration of these two modeling paradigms allows for a more comprehensive understanding of the action of hypoxia-activated drugs. The 3D model can provide spatially resolved oxygen concentrations that drive drug activation within different tumor regions.
The effects of cell death, potentially simulated using a surrogate approach coupled with the local drug concentration, can then be reflected in terms of oxygen consumption within the 3D model, capturing the dynamic interaction between drug activation and the evolving oxygen landscape.
This hybrid strategy aims to balance the need for spatial accuracy in representing the TME with the computational tractability required for simulating complex biochemical reactions and cellular responses. Moreover, the rationale for employing a multiscale modeling framework lies in the trade-off between capturing essential spatial heterogeneities and maintaining computational tractability. 
A pure surrogate (0D) model would neglect spatial gradients of oxygen and drug, which are fundamental to the activation of hypoxia-activated prodrugs. 
However, a fully resolved 3D--1D description of the entire system would be computationally demanding and would complicate parameter estimation, making systematic studies unfeasible. 
The proposed surrogate formulation combines the advantages of both approaches: the 0D component efficiently describes systemic pharmacokinetics, while the 3D-1D module resolves the tumor microenvironment where the oxygen-drug interaction is most critical. 
This framework allows us to address specific scientific questions that neither model could capture alone, such as the influence of intratumoral oxygen heterogeneity on drug activation, the coupling between systemic pharmacokinetics and local tissue hypoxia, and the identification of treatment scenarios where local versus systemic effects dominate.

Therefore, the main purpose of this work is to develop such a multiscale mathematical model capable of capturing the intricate interplay between the transport and activation of hypoxia-activated drugs and the dynamic changes in oxygen availability driven by cell survival. 
Importantly, this multiscale modeling framework enables the application of global sensitivity analysis methods to unravel the influence of key model parameters, including drug properties, physiological factors of the TME, and microvascular characteristics, on drug distribution, activation, and therapeutic efficacy.
Our central hypothesis is that explicitly modeling the pharmacokinetics, transvascular exchange, and hypoxia-dependent activation of tirapazamine within a multiscale hemodynamic framework can provide quantitatively accurate predictions of drug distribution and cytotoxic effect across heterogeneous tumor regions. We further hypothesize that this mechanistic approach offers a foundation for future model generalizations, including fractional-order formulations that may capture additional memory effects in drug–tissue interactions.
Using this sensitivity analysis, we aim to identify the model parameters that most significantly affect relevant outputs, such as the cell survival fraction. Ultimately, this study seeks to provide valuable information for designing and optimizing hypoxia-targeted cancer therapies by identifying the factors that have the greatest impact on treatment outcomes.

\section{Multiscale Model of Hypoxia-Activated Drug Pharmacokinetics and Pharmacodynamics}
\label{mvasMod}

This section details the interconnected models essential for capturing the pharmacokinetics (PK) and pharmacodynamics (PD) of hypoxia-activated drugs, using Tirapazamine (TPZ) as a representative example. Specifically, our surrogate framework integrates: (i) a 3D-1D model of blood flow in the microvascular network coupled with interstitial fluid dynamics, which addresses the heterogeneity of the tumor microenvironment (TME); (ii) a 3D-1D model of oxygen transport, diffusion, and metabolization, crucial for understanding the oxygen dependency of drug activation; and (iii) a model of drug delivery and metabolization within both the vascular and tissue compartments, accounting for its interaction with the oxygen landscape. 
Sections~\ref{sec:flowmod} and~\ref{OxyMod} summarize key modeling components previously developed and validated in authors' previous works~\cite{Possenti2019, Possenti2019b, Possenti2024}, including the coupled 3D-1D hemodynamic model and the reaction–diffusion equations for oxygen transport. These elements are reported here in concise form to provide a self-contained presentation and to enable readers to reproduce the complete computational framework. The present study extends these prior works in important directions: (i) incorporation of a mechanistic model for the pharmacokinetics and tissue distribution of tirapazamine, including plasma clearance and transvascular exchange represented via a semi-permeable membrane model; (ii) introduction of a reaction term describing the bioreductive activation of tirapazamine under hypoxic conditions, which couples drug dynamics to the local oxygen partial pressure; and (iii) formulation of a surrogate (0D) model that reduces the computational cost of simulating drug and oxygen dynamics while preserving the dominant transport and reaction mechanisms. Taken together, these contributions provide a novel, fully coupled 3D-1D-0D \textit{in silico} platform to investigate the efficacy of hypoxia-activated prodrugs in realistic, patient-specific vascular geometries.
The following subsections will detail the mathematical formulation of these interconnected components, focusing on the necessary boundary conditions, parameter considerations, and the physiological interplay that governs drug behavior and efficacy in the hypoxic conditions prevalent in solid tumors. 
This comprehensive model aims to capture the intricate feedback between drug activation and oxygen availability driven by cell survival, as discussed in the Introduction. Appendix~\ref{app-A} and Appendix~\ref{app-B} provide a complete nomenclature, where all model parameters are summarized together with their definition and biological meaning, as well as the acronyms and symbols used throughout the manuscript.

\subsection{A spatially distributed 3D-1D model of the vascular microenvironment}

In the proposed model, the domain $\Omega$ represents a portion of biological tissue (submillimeter) composed of two regions ($\Omega = \Omega_t \bigcup \Omega_v$): the tissue interstitium $\Omega_t$ and the microvascular bed $\Omega_v$. 
$\Omega_t$ is a porous medium, while $\Omega_v$ is an oriented network composed of a set of $N$ cylindrical channels. 
This network is endowed with three sets of variables that indicate the outer surface $\Gamma = \{ \Gamma_i, i = 1,...,N\}$; the radius $R = \{R_i, i = 1,...,N\}$ and the position of the centerline along with the orientation ${\bf \Lambda} = \{{\bf \Lambda_i}, i = 1,..., N\}$ of the selected channel. 
As such, on the vascular bed, the arc length coordinate $s$ is defined as increasing accordingly to the orientation of ${\bf \Lambda_i}, i = 1,..., N$. 
The boundary conditions that complement the problems are imposed at the inlets and outlets, respectively, $\partial\Lambda_{in}$ and $\partial\Lambda_{out}$. 
Since we approximated the vascular domain to a 1D domain, from now on, the microvascular domain refers to ${\bf \Lambda}$ and the tissue domain by \(\Omega\), with \(\Omega \simeq \Omega_t\).
In what follows, for the compact description of the governing equations of flow and transport in the vascular microenvironment, we will denote the geometrical data, defining the domains \(\Omega\) and \(\Lambda\) as \(\mathcal{D}\).

\subsubsection{The Microvascular Flow Model}
\label{sec:flowmod}

Blood flow is crucial to understanding how drugs are distributed within the vascular network and subsequently delivered to tumor tissues. 
The rate and pattern of blood flow determine the delivery of oxygen and the drug itself to the hypoxic regions. 
The blood flow can be modeled using fluid dynamics principles, often employing equations like Poiseuille's law for vascular flow and Darcy's law for tissue perfusion. 
This involves understanding how blood pressure, vascular resistance, and tissue permeability affect drug delivery.
The mathematical model describing flow dynamics and hematocrit transport in a vascular network is represented here by the combined framework \(\mathcal{F} \& \mathcal{H}\). 
This framework provides a comprehensive, yet simplified, description of fluid dynamics and red blood cell distribution within the system:
\[
\mathcal{F} \& \mathcal{H}(p, \mathbf{u}, H; \, \mathcal{D}, P_v^{in}, \theta) = 0\, .
\]
The notation in the equation above distinguishes between the unknowns and the parameters of the problem. 
In this expression, the variables \(p\), \(\mathbf{u}\), and \(H\) before the semicolon (;) represent the unknowns. 
The terms after the semicolon —\(\mathcal{D}\), \(P_v^{in}\), and \(\theta\)— represent the parameters of the problem. 
\(\mathcal{D}\) includes domain characteristics, \(P_v^{in}\) is the parameter of the boundary condition on the input of the vascular network, and \(\theta\)  represents any additional parameters influencing the system.
This notation separates the variables to be solved from the model's fixed parameters and will be consistently used throughout the document for all other abstract models.
The fluid dynamics model, \(\mathcal{F}\), incorporates both the Poiseuille flow within the vasculature and the Darcy flow in the tissue. 
The Poiseuille flow within the vasculature and the Darcy flow in the tissue are nonlinearly coupled via extravasation and lymphatic drainage terms. 
The continuity equations govern both domains (\(\Lambda\) and \(\Omega\)).
The complete model is expressed as:
\begin{equation} \label{eq:F}
\mathcal{F}:=
\begin{cases}
\nabla\cdot{\bf u}_t+L_{p}^{L F}{\frac{S}{V}}(p_{t}-p_{L})-f_b(\nobar{p}_t,p_v)\,\delta_{\Lambda} & \mathrm{in}\quad\Omega\\
\mathbf{u}_{t}+{\frac{\kappa}{\mu_t}}\nabla p_{t} & \mathrm{in}\quad \Omega\\
\partial_s \left( \pi R^2 u_v \right)+f_b(\nobar{p}_t,p_v) & \mathrm{in}\quad \Lambda\\
8\mu_v u_v+R^{2} \partial_s p_{v} & \mathrm{in}\quad \Lambda\\
p_v-(p_0+\Delta p) & \mathrm{in}\quad \partial\Lambda_{in}\\
p_v-p_0 & \mathrm{in}\quad \partial\Lambda_{out}\\
{\bf u}_t\cdot\boldsymbol{n} & \mathrm{in}\quad \partial\Omega\\
\end{cases}
\end{equation}
In this formulation, the quantities with subscripts \(t\) and \(v\) refer to the tissue and the vascular bed, respectively. 
Here, \({\bf u}_t\) and \(u_v\) are fluid velocities, while \(p_t\) and \(p_v\) are pressures. 
The pressures \(p_0\) and \(\Delta p\) correspond to the outlet pressure and the pressure difference between inlets and outlets (therefore, we define the inlet pressure as \(p_0 + \Delta p\)).
The viscosities \(\mu_t\) and \(\mu_v\) denote the dynamic viscosities of the two fluids.
Furthermore, the term \(L_{p}^{L F}\frac{S}{V}(p_{t}-p_{L})\) represents the volumetric flow rate due to lymphatic drainage, with \(L_{p}^{L_F}\) being the hydraulic permeability of the lymphatic walls.
The function \(f_b(\nobar{p}_t,p_v)\) models fluid extravasation according to the Starling model:
\begin{equation}
f_b(\nobar{p}_t,p_v)=2\pi R L_p[(p_v-\nobar{p}_t)-\sigma(\pi_v-\pi_t)]\, ,
\end{equation}
where \(L_p\) is the hydraulic conductivity, \(\pi_v\) and \(\pi_t\) are the osmotic pressure gradients across the capillaries, and \(\sigma\) is the reflection coefficient.
Model \(\mathcal{F}\) is further extended with the coupled one-dimensional red blood cell (RBC) transport model \(\mathcal{H}\).
The hematocrit transport model, \(\mathcal{H}\), ensures the conservation and distribution of RBC concentration (i.e., hematocrit, \(H\)) throughout the vascular network, maintaining mass balance. 
The governing equations of the hematocrit \(H\) within \(\Lambda\) read as follows:
\begin{equation}
    \label{eq:1D_H}
\mathcal{H}:=
\begin{cases}
\pi R^2 u_v \partial_s H  - f_b(\nobar{p}_t,p_v) H  & \mathrm{in} \quad \Lambda\\
H - H_{in} & \mathrm{in}\quad \partial\Lambda_{in}\\
\partial_s H & \mathrm{in}\quad \partial\Lambda_{out}\\
\end{cases}
\end{equation}
Here, \(H_{in}\) is the input hematocrit value. 
This model assumes hematocrit as a conserved quantity, which means that RBCs do not extravasate from \(\Lambda\) and are not degraded during transport.
Furthermore, network connectivity allows only for bifurcations or anastomoses, ensuring mass conservation at all junctions, leveraging the problem closure defined by Pries et al.~\cite{Pries2005, Possenti2019}.
Finally, we note that the microvascular flow model focuses on equilibrium conditions of spatial distributions, ignoring any time-dependent dynamics that might occur in the system. 
This simplification is useful for studying long-term behavior and the overall distribution of flow and hematocrit in the vascular network, but it can overlook transient events and fluctuations that could be important in certain physiological or pathological situations. 

\subsubsection{Oxygen Transport and Metabolism}
\label{OxyMod}

The efficacy of hypoxia-activated prodrugs is directly related to the level of oxygen in the tissue, since hypoxic conditions activate these drugs.
Oxygen transport models are used to simulate the diffusion and consumption of oxygen within the tissue by solving diffusion equations that consider tissue oxygen demand and capillary oxygen supply. 
Note that the interaction between the flow and the oxygen model is characterized by a one-way coupling because \(\mathcal{F} \& \mathcal{H}\) influences oxygen transport, but not vice versa.
The oxygen transport model describes the diffusion and metabolism of oxygen within tissue and vasculature, characterized by oxygen concentration \( c^{ox} \):
\[
\partial_t c_t^{ox} 
+ \mathcal{T}_{ox}(c_t^{ox}; \mathcal{D}, p, \mathbf{u}, H, c_t^{ox,in}, \theta^{ox}) 
+ \mathcal{M}_{ox}(c_t^{ox})
= 0\, ,
\]
where \( c^{ox, in} \) represents the boundary condition for oxygen concentration at inflow and \(\theta^{ox}\) includes all physical parameters, such as viscosities, hydraulic conductivities, and reflection coefficients.
To prepare the coupling with TPZ, we have the effects related to oxygen transport, represented by \(\mathcal{T}_{ox}\), from those related to oxygen metabolization, denoted by \(\mathcal{M}_{ox}\). 
This compact notation hides several phenomena that govern oxygen transport and delivery to cells, both in the \(\Lambda\) and the \(\Omega\) domains. 
In the vasculature, the concentration of oxygen available in the blood is the sum of the concentration of dissolved oxygen $c^{ox}_v$ and hemoglobin-bound oxygen $c^{ox}_{HbO_2}$. 
Moreover, given the fast kinetics, we neglect transient phenomena related to hemoglobin binding, so that $c^{ox}_v$ and $c^{ox}_{HbO_2}$ are always chemically balanced.
Consequently, $c^{ox}_{HbO_2}$ is a function of $c^{ox}_v$:
\begin{equation}
\label{eq:hill_C}
{c^{ox}_{HbO_2}(c^{ox}_v) = k_1 H \frac{c^{ox, \gamma}_v}{c^{ox, \gamma}_v + (\alpha_{pl}\, p_{s_{50}})^\gamma}}\, ,
\end{equation}
with $k_1$ a constant given by the H\"ufner factor $N$ times the Mean Corpuscular Hematocrit Concentration $MCHC$; $\alpha_{pl}$ the solubility of oxygen in plasma; $p_{s_{50}}$ the oxygen partial pressure at hemoglobin half-saturation; and $\gamma$ the Hill exponent. 
The evolution of the oxygen concentration in $\Lambda$ is given by:
\begin{equation}
-\pi R^2 D_v^{ox} \partial^2_s c^{ox}_v 
+\pi R^2 \partial_s \left(u_v c^{ox}_v + u_v k_1 H \frac{c^{ox, \gamma}_v}{c^{ox, \gamma}_v + k_2}\right)  = -  f_c^{ox}(\nobar{p}_t,p_v,\nobar{c}^{ox}_t,c^{ox}_v)\quad \textrm{on $\Lambda$} \ ,
\label{eq:strong_Ov}
\end{equation}
where $D_v^{ox}$ is oxygen diffusion coefficient in $\Lambda$; $k_2 = (\alpha_{pl} p_{s_{50}})^\gamma$ and $f_c^{ox}(\nobar{p}_t,p_v,\nobar{c}^{ox}_t,c^{ox}_v)$ is a coupling term modeling the diffusion of $c^{ox}_v$ from $\Lambda$ to $\Omega$. 
In this work, we adopt the Kedem-Katchalsky model:
\begin{equation}
\label{eq:KK}
f_c^{ox}(\nobar{p}_t,p_v,\nobar{c}^{ox}_t,c^{ox}_v)  = 2 \pi R P^{ox} (c^{ox}_v - \nobar{c}^{ox}_t) + (1-\sigma^{ox})~\left(\frac{c^{ox}_v + \nobar{c}^{ox}_t}{2}\right)\cdot f_b(\nobar{p}_t,p_v)\, ,
\end{equation}
with $c^{ox}_t$ the tissue oxygen concentration and $\nobar{c}^{ox}_t$ its mean value on $\Omega$; $P^{ox}$ the permeability of the vascular wall to oxygen, $\sigma^{ox}$ a reflection coefficient relative to the oxygen molecule~\cite{jarzynska2006application}.
Although this modeling approach describes both diffusive and advective oxygen flow, we remark that oxygen delivery from microvasculature is predominantly a diffusion-dominated problem ($ 2 \pi R P^{ox} (c^{ox}_v - \nobar{c}^{ox}_t) \gg (1-\sigma^{ox})~\left(\frac{c^{ox}_v + \nobar{c}^{ox}_t}{2}\right)\cdot f_b(\nobar{p}_t,p_v)$).
Note that the prescribed evolution for $c^{ox}_v$ in Eq.~\eqref{eq:strong_Ov} holds by assuming $c^{ox}_v$ as conserved in $\Lambda$ (no oxygen consumption within blood flow) and that the diffusion coefficient for $c^{ox}_{HbO_2}$ is null.
On the other hand, transport and diffusion of oxygen concentration in $\Omega$ corresponds to:
\begin{equation}
- \nabla \cdot \left(c^{ox}_t\mathbf{u}_t - D_t^{ox} \nabla c^{ox}_t \right) 
= f_c^{ox}(\nobar{p}_t,p_v,\nobar{c}^{ox}_t,c^{ox}_v) \delta_\Lambda - m(c^{ox}_t)\quad \textrm{on $\Omega$}\, ,
\label{eq:general_Ot}
\end{equation}
where $D_t^{ox}$ is oxygen diffusion coefficient in $\Omega$ and $\mathbf{u}_t$ the fluid velocity.
$m(c^{ox}_t)$ is the rate of oxygen depletion due to the metabolic activity of the tissue (Michaelis-Menten model):
\begin{equation}
\label{eq:MichaelisMentenConc}
{m(c^{ox}_t) = V_{max} \frac{c^{ox}_{t}}{c^{ox}_{t} + \alpha^{ox}_t p_{m_{50}}}}\, .
\end{equation}
with $V_{max}$ its maximum consumption rate; $p_{m_{50}}$ its partial pressure at half consumption rate; and $\alpha^{ox}_t$ its solubility in the tissue. As a result, in the conditions where the metabolism of oxygen is not affected by TPZ, the term \(\mathcal{M}_{ox}\) is given by:
\[
\mathcal{M}_{ox}(c_t^{ox}; \theta^{ox})={m(c^{ox}_t) = V_{max} \frac{c^{ox}_{t}}{c^{ox}_{t} + \alpha^{ox}_t p_{m_{50}}}}\,.
\]
The overall proposed model for oxygen transport and diffusion is:
\begin{equation} \label{eq:O}
\mathcal{T}_{ox}:
\begin{cases}
- \nabla \cdot \left(c^{ox}_t\mathbf{u}_t - D_t^{ox} \nabla c^{ox}_t \right) 
- f_c^{ox}(\nobar{p}_t,p_v,\nobar{c}^{ox}_t,c^{ox}_v)\delta_\Lambda & \mathrm{in}\quad\Omega\\
-\pi R^2 D_v^{ox} \partial^2_s c^{ox}_v 
+\pi R^2 \partial_s \left(u_v c^{ox}_v + u_v k_1 H \frac{c^{ox, \gamma}_v}{c^{ox, \gamma}_v + k_2}\right)  
+  f_c^{ox}(\nobar{p}_t,p_v,\nobar{c}^{ox}_t,c^{ox}_v) & \mathrm{in}\quad \Lambda\\
c^{ox}_v - c^{ox}_{in} & \mathrm{in}\quad \partial\Lambda_{in}\\
-D_v^{ox}\partial_s c^{ox}_v & \mathrm{in}\quad \partial\Lambda_{out}\\
-D_t^{ox}\nabla c^{ox}_t\cdot\boldsymbol{n}-\beta_{ox} (c^{ox}_t-c^{ox}_0)& \mathrm{in}\quad \partial\Omega\\
\end{cases}
\end{equation} 
At $\partial\Lambda_{in}$ the oxygen concentration $c^{ox}_{in}$ is specified. For the tissue, we simulate the presence of an adjacent tissue domain with boundary conductivity $\beta_{ox}$ and far-field concentration $c^{ox}_0$.
In fact, the latter is only one-way coupled with $\mathcal{F}$ and $\mathcal{H}$ through $u_v$, ${\bf u}_t$, and $H$, while $c^{ox}_v$ and $c^{ox}_t$ have no influence on the blood dynamics.
We note that when oxygen metabolism is not affected by hypoxia-activated drugs, such as TPZ, the oxygen level may reach a steady state, determined by the equilibrium between the supply of microvessels and the consumption of cells~\cite{Possenti2021}. Under these conditions, the oxygen model is steady and can be simplified as follows:
\begin{equation}\label{eq:Ofull}
    \underbrace{\mathcal{T}_{ox}(c^{ox}; \, \mathcal{D}, p, \mathbf{u}, H, c^{ox}_{in}, \theta^{ox}) 
+ \mathcal{M}_{ox}(c^{ox}; \theta^{ox})}_{\mathcal{O}}
= 0\, .
\end{equation}

\subsubsection{Pharmacokinetics and pharmacodynamics of Tirapazamine}
\label{phakinMod}

Accurate drug delivery models are essential to predict the concentration of drugs reaching the target tissue, directly affecting their therapeutic efficacy and safety profile. 
The activation of hypoxia-activated drugs like Tirapazamine (TPZ) depends on their metabolic conversion, considering variations in metabolic rates under different oxygen levels. 
Understanding these kinetics is crucial to predicting drug activation and efficacy.
TPZ is a prototypical hypoxia-activated prodrug whose pharmacokinetics and pharmacodynamics have been extensively characterized. Early preclinical studies demonstrated that TPZ exhibits selective activation under hypoxic conditions, with tumor-to-plasma ratios consistent with its hypoxia-dependent cytotoxicity~\cite{Dorie1984, Brown1993}. Clinical pharmacokinetics investigations confirmed this tissue selectivity and provided evidence for heterogeneous metabolism across tumor regions~\cite{Patterson2000}. More recent work has highlighted the limited penetration of TPZ into poorly perfused tumor zones, an aspect that critically affects its therapeutic efficacy~\cite{PKPD}. Pharmacokinetics/pharmacodynamics modeling studies further established the role of spatial distribution and hypoxia gradients in shaping drug activity, and identified strategies to improve tissue penetration and cytotoxic potency in TPZ analogues~\cite{Hicks2010}. Together, these studies provide the biological basis for developing the following model of TPZ concentration in hypoxic tumor tissue.

The TPZ pharmacokinetic and pharmacodynamic model describes the distribution and effects of the tirapazamine drug within the vascular network and tissue:
\begin{equation}
\partial_t c^{tpz} + \mathcal{TPZ}(c^{tpz}; 
\mathcal{D}, p, \mathbf{u}, H,  SF, c^{ox}, 
c_{v,in}^{tpz}, \theta^{tpz})
=0\, ,
\label{mod:TPZ}
\end{equation}
where \(c^{tpz}=[c_{v}^{tpz},c_{t}^{tpz}]\) is the concentration of tirapazamine in the vascular bed and the tissue, respectively. \(c_{v, in}^{tpz}\) denotes the drug concentration at the inflow boundary, and \(\theta^{tpz}\) includes drug-specific parameters, such as diffusion coefficients and metabolic rates. 
When tirapazamine is modeled within the tumor microenvironment, the dynamic nature of the cellular population of the tumor must be taken into account. Consequently, the TPZ pharmacokinetic and pharmacodynamic model must incorporate temporal changes in cellular density and distribution to accurately simulate drug activation and efficacy. 
This model must dynamically link changes in cell population with fluctuations in oxygen availability and TPZ activation, ensuring that they reflect the complex and evolving nature of the tumor microenvironment. 
This is done by introducing a new variable in the model, named \(SF\), that is the surviving fraction of cells that metabolize tirapazamine. Let us subdivide the model $\mathcal{TPZ}$ into two parts, corresponding to the drug transport and the drug metabolization by cells. 
As such, we rewrite Eq.~\eqref{mod:TPZ} as follows:
\begin{equation}
\partial_t c^{tpz} + \mathcal{T}_{tpz}(c^{tpz}; \mathcal{D}, p, \mathbf{u}, H,
c_{v,in}^{tpz}, \theta^{tpz})
+ \mathcal{M}_{tpz}(c_{t}^{tpz}, SF; c_t^{ox})
=0\, .
\label{eq:pharmacokinetics}
\end{equation}
Also, note that tirapazamine is a hypoxia-activated drug whose metabolism depends on the local oxygen concentration.
The operator $\mathcal{TPZ}$ involves both the vasculature and tissue domains.
The concentration of injected drug $c_v^{tpz}$ evolves in the vascular bed governed by the following boundary value problem:
\begin{equation}
\begin{cases}
\partial_t c_{v}^{tpz} + \partial_s \left(c_v^{tpz} u_v
-D_v^{tpz} \partial_s c_v^{tpz}\right)=
-\frac{1}{\pi R^2}f_c^{tpz}(\nobar{p}_t,p_v,\nobar{c}_t^{tpz},c_v^{tpz}) & \mathrm{in}\quad \Lambda \\
c_v^{tpz}=c_{v,in}^{tpz} & \mathrm{in}\quad \partial\Lambda_{in}\\
-D_v^{tpz} \partial_s c_v^{tpz}=0 & \mathrm{in}\quad \partial\Lambda_{out}\\
\end{cases}
\end{equation}
with $t$ time coordinate and $D_v^{tpz}$ the diffusion coefficient in $\Lambda$ and $c_{v, in}^{tpz}$ is a function returning the injected dose of drug concentration in inlet $\partial \Lambda_{in}$,
\begin{equation}
c_{v,in}^{tpz} = \begin{cases}
c_v^{tpz} = a\, t & \mathrm{for}\quad t\in(0,T_{P1})\\
c_v^{tpz}=c_{v_0}^{tpz} & \mathrm{for}\quad t\in(T_{P1},T_{P2})\\
c_v^{tpz}=c_{v_0}^{tpz}(e^{-\frac{CL}{V}\cdot t})& \mathrm{for}\quad t\in(T_{P2},T_{P2}+5\cdot\frac{V}{CL})\\
c_v^{tpz}=0 & \mathrm{for}\quad t>T_{P2}+5\cdot\frac{V}{CL}\\
\end{cases}
\end{equation}
{where $T_{P1}$ corresponds to the time at which $c_{v, in}^{tpz}$ reaches a plateau, $c_{v0}^{\mathrm{tpz}}$ is the TPZ concentration value in the vascular compartment at this plateau (i.e. an effective delivered concentration), while $T_{P2}$ denotes the duration of the plateau, before the effective concentration decays with exponential rate.}
The increase in vascular concentration is assumed to be linear, with a slope of $a$. 
We have accounted for a decay time of $5\tau = 5\frac{V}{CL}$, {being $V$ the apparent volume of distribution, $CL$ the systemic clearance rate, and $\tau = V / CL$ the characteristic decay time of TPZ in plasma.} 
 On the tissue side, we consider the following model~\cite{PKPD}:
\begin{equation}
\begin{cases}
\partial_t c_{t}^{tpz} 
+ \nabla\cdot (c_t^{tpz}\boldsymbol{u}_t - D_t^{tpz}\nabla c_t^{tpz}) + L_{p}^{LF} \frac{S}{V}(p_t-p_L)c_t^{tpz} = f_c^{tpz}(\nobar{p}_t,p_v,\nobar{c}_t^{tpz},c_v^{tpz})\delta_{\Lambda}+\phi m^{tpz}(c_t^{tpz}) &\mathrm{in} \quad \Omega\\
-D_t^{tpz}\nabla c_t^{tpz}\cdot\boldsymbol{n}=\beta^{tpz} (c_t^{tpz}-c_0^{tpz})& \mathrm{in}\quad \partial\Omega\\
\end{cases}
\end{equation}
where $c_t^{tpz}$ is the drug concentration in $\Omega$ {and $f_{c}^{\text{tpz}}$ represents the mass flux of TPZ across the capillary wall, modeled analogously to oxygen exchange in Eq.~\eqref{eq:KK} by a semi-permeable membrane model.}
In the same fashion as per $\mathcal{O}$, the problem is complemented by the conductive boundary condition ($\beta^{tpz}$ being the conductivity of the walls and the concentration of the far field $c_0^{tpz}$).
To model the diffusivity of TPZ in tissue, we rely on empirical expressions derived from \textit{in vitro} studies on multicellular cancer layers (MCLs)~\cite{zou2020application}, which aim to predict drug diffusivity based on molecular descriptors. In this framework, the diffusion coefficient \( D_t^{tpz} \) used in our model is estimated from the following relation:
\begin{equation}
    \log(D_t^{tpz}) = a + b \log(MW) + \frac{c}{1 + \exp\left( \frac{\log P_{7.4} - x + y \cdot HD + z \cdot HA}{w} \right)}\, ,
    \label{eq:D}
\end{equation}
where \( MW \) is the molecular weight, \( \log P_{7.4} \) is the octanol/water partition coefficient at pH 7.4, \( HD \) and \( HA \) denote the numbers of hydrogen bond donors and acceptors, respectively. The coefficients \( a, b, c, w, x, y, z \) are empirical parameters fitted to the experimental data and capture how physicochemical properties influence diffusion in multicellular layers. For this reason, the estimated value of \( D_t^{tpz} \) can be interpreted as an effective diffusion coefficient in multicellular layer (MCL) in-vitro tissue environments, often denoted in the literature as \( D_{\text{mcl}} \).
For assessing TPZ metabolization, we consider the combination of two terms, precisely \(\phi\) and \( m^{tpz}(c_t^{tpz})\), where \(\phi\) represents the population of viable cells and \(m^{tpz}(c_t^{tpz})\) is the rate of drug metabolism. 
We rewrite the viable cell population introducing the surviving fraction \(SF\), defining $\displaystyle \phi=\phi_0 SF$, with $\phi_0$ corresponding to the initial cellular volume fraction and $SF$ to the cell surviving fraction under the action of the drug.
The term $m^{tpz}(c_t^{tpz})$ is then defined by a modified Michaelis-Menten dynamics with an effective term depending on oxygen concentration to describe the hypoxia-activated drug behavior:
\begin{equation}
\label{metcttpz}
m^{tpz}(c_t^{tpz},c^{ox}) = \left(\frac{K}{K+c^{ox}_t}\right)\left(k_{met} \, {c_{t}^{tpz}}{+}\ \frac{V^{tpz}_{max} \,{c_{t}^{tpz}}}{K_m\ {+}\ {c_{t}^{tpz}}}\right)\, , 
\end{equation}
where $k_{met}$ is the first order metabolic rate constant, $V^{tpz}_{max}$ is the maximal rate of Michaelis-Menten metabolism, $K_m$ is the Michaelis constant and $K$ represents the oxygen concentration to halve $m^{tpz}(c^{tpz}_t)$~\cite{pruijn2005tpz}. 
As a result, the TPZ metabolization model becomes:
\begin{equation}
\mathcal{M}_{tpz}(SF, c^{ox}, c_{t}^{tpz}) =\phi_0\, SF \, m^{tpz}(c_t^{tpz}, c_{t}^{ox}) \, ,
\label{mod:MTPZ}
\end{equation}
$\phi_0$ being the initial cellular volume fraction.
To describe the surviving fraction $SF$, we here introduce the pharmacodynamics model defining the drug's effect on cancer cells.
Due to the tirapazamine's action, the population of viable cells in the system is not constant. 
For this reason, the surviving fraction of cancer cells is regulated by an exponential law. More precisely, the rate of the logarithm of $SF$ is modeled as a linear function of $m^{tpz}$ and $c_t^{tpz}$~\cite{PKPD,siim1996tirapazamine}:
\begin{equation}
    {-}\ \frac{d\ \mathrm{log}\ SF}{dt} = \alpha\, c_{t}^{tpz} \, m^{tpz}(c_t^{tpz}, c_{t}^{ox})\, .
    \label{dlogSF}
\end{equation}
$\alpha$ is a constant heuristically derived from linear regression of experimental data quantifying the cytotoxic potency of TPZ.
As a consequence:
\begin{equation}
    SF(c_t^{tpz}, c_{t}^{ox}) = \exp{\left(- \int_0^t \alpha\, c_{t}^{tpz}\ m^{tpz}(c_t^{tpz}, c_{t}^{ox}) \ d\tau\right)}\, .
    \label{defSF}
\end{equation}

In general, we formulate a model to describe both pharmacokinetics and pharmacodynamics. 
The pharmacokinetics is described by equation \ref{eq:pharmacokinetics}, where $\mathcal{M}_{tpz}(SF, c^{ox}, c_{t}^{tpz})$ is defined in {Eq.~\eqref{mod:MTPZ}}, and $\mathcal{T}_{tpz}$ is:
\begin{equation}
\label{TPZ}
\mathcal{T}_{tpz}:
\begin{cases}
\nabla\cdot (c_t^{tpz}\boldsymbol{u}_t - D^{tpz}_t\nabla c^{tpz}_t) 
+ L_{p}^{LF} \frac{S}{V}(p_t-p_L)c^{tpz}_t
- f_c^{tpz}(\nobar{p}_t,p_v,\nobar{c}^{tpz}_t,c^{tpz}_v)\delta_{\Lambda} & \mathrm{in}\quad\Omega\\
\partial_s \left(c^{tpz}_v u_v -D_v^{tpz} \partial_s c^{tpz}_v\right)
+\frac{1}{\pi R^2}f_c^{tpz}(\nobar{p}_t,p_v,\nobar{c}^{tpz}_t,c^{tpz}_v) & \mathrm{in}\quad \Lambda\\
p_v-p_0-\Delta p & \mathrm{in}\quad \partial\Lambda_{in}\\
p_v-p_0 & \mathrm{in}\quad \partial\Lambda_{out}\\
c_v^{tpz} - c_{v,in}^{tpz} & \mathrm{in}\quad \partial\Lambda_{in}\\
\partial_s c^{tpz}_v& \mathrm{in}\quad \partial\Lambda_{out}\\
-D_t^{tpz}\nabla c_t\cdot\boldsymbol{n}-\beta (c^{tpz}_t-c^{tpz}_0)& \mathrm{in}\quad \partial\Omega
\end{cases}
\end{equation}

The effect of the drug on cell viability also influences the uptake of oxygen, altering the oxygen concentration and, consequently, affecting the drug activity. Consistent with the modeling approach adopted for TPZ, we similarly define the oxygen metabolism term as:
\begin{equation}
\mathcal{M}_{ox}(SF, c^{ox}) = \phi_0 \, SF \, m^{ox}(c^{ox}),
\end{equation}
where the surviving fraction $SF$ is calculated according to the tissue TPZ concentration $c_{t}^{tpz}$ as detailed in equation~\eqref{defSF}. Thus, the modified oxygen transport and metabolism equation becomes:
\begin{equation}
\partial_t c^{ox} + \mathcal{T}_{ox}(c^{ox}; \mathcal{D}, p, \mathbf{u}, H, c^{ox}_{in}, \theta^{ox}) + \phi_0 \, SF \, m^{ox}(c^{ox})\,.
\end{equation}

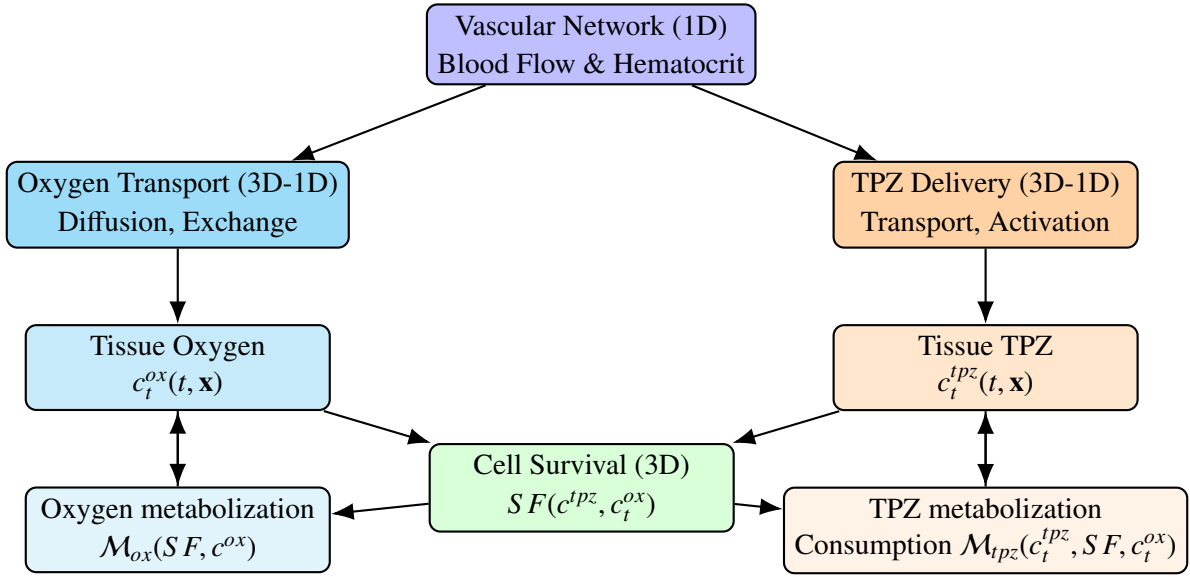
\begin{figure}[ht]
\centering
\begin{tikzpicture}[
  every node/.style={font=\small},
  box/.style={draw, rounded corners, thick, align=center, minimum width=4.0cm, minimum height=1.0cm},
  arrow/.style={thick, -{Latex[length=3mm]}},
  doublearrow/.style={thick, <->, >=Latex},
  >=Latex,
  node distance=1.0cm and 1.0cm
]

% Top layer: Vascular network
\node[box, fill=blue!25] (vascular) {Vascular Network (1D)\\[1pt] Blood Flow \& Hematocrit};

% Middle left: Oxygen transport
\node[box, below left=of vascular, fill=cyan!35] (oxygen) {Oxygen Transport (3D-1D)\\[1pt] Diffusion, Exchange};

% Middle right: Drug delivery
\node[box, below right=of vascular, fill=orange!35] (drug) {TPZ Delivery (3D-1D)\\[1pt] Transport, Activation};

% Lower left: Oxygen in tissue
\node[box, below=of oxygen, fill=cyan!20] (ctox) {Tissue Oxygen\\$c_t^{ox}(t,\mathbf{x})$};

% Lower right: Drug in tissue
\node[box, below=of drug, fill=orange!20] (ctpz) {Tissue TPZ\\$c_t^{tpz}(t,\mathbf{x})$};

% Bottom center: SF
\node[box, below=of $(ctox)!0.5!(ctpz)$, fill=green!15] (sf) {Cell Survival (3D)\\[1pt] $SF(c^{tpz}, c_{t}^{ox})$};

% Bottom center: mox
\node[box, below=of ctox, fill=cyan!10] (mox) {Oxygen metabolization\\[1pt] $\mathcal{M}_{ox}(SF, c^{ox})$};

% Bottom center: mtpx
\node[box, below=of ctpz, fill=orange!10] (mtpz) {TPZ metabolization\\Consumption $\mathcal{M}_{tpz}(c_{t}^{tpz}, SF, c_t^{ox})$};

% Arrows between boxes
\draw[arrow] (vascular) -- (oxygen);
\draw[arrow] (vascular) -- (drug);
\draw[arrow] (oxygen) -- (ctox);
\draw[arrow] (drug) -- (ctpz);
\draw[arrow] (ctox) -- (sf);
\draw[arrow] (ctpz) -- (sf);
\draw[arrow] (ctox) -- (mox);
\draw[arrow] (ctpz) -- (mtpz);
\draw[arrow] (sf) -- (mox);
\draw[arrow] (sf) -- (mtpz);
\draw[arrow] (mox) -- (ctox);
\draw[arrow] (mtpz) -- (ctpz);

\end{tikzpicture}
\caption{Schematic of the full 3D-1D model. The vascular network (1D) supplies oxygen and drugs to the tissue (3D). Spatially resolved concentrations influence the cell survival model at each point in space and time (3D), which feeds back to modulate oxygen and TPZ consumption and reshape the microenvironment.}
\label{fig:hybrid_model_schematic}
\end{figure}

In conclusion, the model we formulate here for the study of hypoxia-activated drugs is the following:
\begin{equation}
\label{modelAll}
    \begin{cases}
     & \mathcal{F} \& \mathcal{H}(p, \mathbf{u}, H; \mathcal{D}, P_v^{in}, \theta) =0
    \\
    & \partial_t c^{ox} 
    + \mathcal{T}_{ox}(c^{ox}; \mathcal{D}, p, \mathbf{u}, H, c^{ox}_{in}, \theta^{ox})
    + \mathcal{M}_{ox}(SF, c^{ox}) = 0
    \\
    & \partial_t c^{tpz} 
    + \mathcal{T}_{tpz}(c^{tpz}; \mathcal{D}, p, \mathbf{u}, H, c_{v,in}^{tpz}, \theta_{T^{tpz}})
    + \mathcal{M}_{tpz}(c_{t}^{tpz}, SF, c_t^{ox}) = 0
    \\
    & SF(c_{t}^{tpz}) = \exp{\left(- \int_0^t \alpha\, c_{t}^{tpz}\, m(c_{t}^{tpz})d\tau\right)}
    \end{cases}
\end{equation}
This model, also illustrated in the schematic of Figure \ref{fig:hybrid_model_schematic}, highlights the essential pharmacokinetics and pharmacodynamics processes, capturing the interaction of tirapazamine with the vascular and tissue environments, influenced by oxygen levels and drug properties. By integrating this model with the fluid dynamics, hematocrit, and oxygen transport models, a comprehensive understanding of drug behavior in hypoxic tumor regions is achieved.

\subsection{A surrogate model for pharmacokinetics and pharmacodynamics: the 0D model}
\label{sec:model0D}

The 3D-1D mixed-dimensional model presented in~Eq.s~\eqref{modelAll} shows a nonlinear interdependence of $SF$, $c_t^{ox}$ and $c_t^{tpz}$, as shown in Figure \ref{fig:hybrid_model_schematic}.
This model definition represents a peculiar feature for properly describing reliable pharmacokinetics, but is incompatible with a sensitivity analysis approach due to computational demands.
To address this, we developed a surrogate model (0D model) incorporating nonlinear dynamics to support the 3D-1D mixed-dimensional model. 
Specifically, our approach is to leverage the linear relationship for drug consumption while thoughtfully incorporating nonlinear dynamics that are \textit{a priori} determined using a suitable 0D model.  
{In the present work, we construct a surrogate model that reduces the complex distributed system governing TPZ pharmacokinetics and pharmacodynamics (described by a set of coupled partial differential equations) to a system of ordinary differential equations formulated in terms of spatially averaged quantities. This approach decreases computational complexity while preserving the dominant transport and reaction dynamics.}
We formulate a surrogate model based on Eq.~\eqref{TPZ} in $\Omega$ to quickly compute the spatial average of $c_t^{tpz}$ and $SF$.

Neglecting the space dependence, the {ordinary differential problem ($d_t$ being the time derivative)} for TPZ reads:
\begin{equation}
d_t c_{t}^{\text{tpz}}+\phi(SF) m^{\text{tpz}}(c_t^{tpz},c_t^{ox}) 
 + L_{p}^{LF} \frac{S}{V}(p_t-p_L)c_t^{\text{tpz}} = f_c^{\text{tpz}}(p_t,p_v,c_t^{\text{tpz}},c_v^{\text{tpz}})\, ,
\end{equation}
where we consider the consumption rate $m^{tpz}(c_t^{tpz})$, the lymphatic drainage ($ L_{p}^{LF} \frac{S}{V}(p_t-p_L)c_t^{tpz}$), and the forcing term {$f_c^{\text{tpz}}$}. 
In this study, the influence of lymphatic drainage and capillary leakage is considered negligible. This assumption is supported by experimental observations that indicate that functional lymphatic vessels are often absent in the core regions of solid tumors, resulting in effective lymphatic clearance rates as low as \(10^{-6}\ \text{s}^{-1}\) or less~\cite{Baxter1989a, Baxter1990a}. Similarly, although capillary leakage is typically enhanced in tumors, the time scales associated with solute extravasation are much longer than the characteristic diffusion and metabolism times considered here. Therefore, their contribution to drug transport dynamics is deemed subdominant and neglected in the reduced-order formulation.
By neglecting lymphatic drainage and capillary leakage, three main contributions are identified:
\begin{itemize}
\item[$i)$] the flux of the drug due to the permeability of the capillary walls: $\frac{S}{V}P(c^{tpz}_v-c_t^{p,tpz})$;
\item[$ii)$] the diffusion flux from the perivascular environment to the tissue: $\frac{S}{V}\frac{D_t^{tpz}}{L}(c_t^{p,tpz}-c^{tpz}_t)$;
\item[$iii)$] drug consumption in the tissue.
\end{itemize}
{Here, $P$ is the permeability of the capillary wall to TPZ and $L$ is the representative distance for the diffusion of the solute between the perivascular space and the bulk of the tissue (characteristic diffusion path length, $2.5 \mu m$). The quantity $S/V$ denotes the vessel surface area per unit tissue volume.} Moreover, we introduce a new variable $c_t^{p,tpz}$ representing the concentration of drugs in the perivascular environment to accurately describe diffusion in the tissue without explicitly including the spatial coordinate.
Consequently, $c_t^{tpz}$ is the average concentration in the tissue and $D_t^{tpz}$ is the solute diffusion coefficient.
Thus, the homogeneous problem holds:
\begin{equation}
     d_t c_t^{tpz} + \phi(SF) m^{tpz}(c_t^{tpz},c_t^{ox}) =  P \frac{S}{V} (c_v^{tpz} - c_t^{p,tpz}) + \frac{D^{tpz}_t}{L} \frac{S}{V} (c_t^{p,tpz}-c_t^{tpz})\, .
 \end{equation}

In the same fashion as in electromagnetism, the flux of the drug from the vessel to the tissue can be modeled as a \textit{current} and the concentration differences as a \textit{voltage}, obtaining the resulting \textit{resistance} exerted by the drug~\cite{Bird2002Transport, PROSI2005903}. 
Let $\phi_1$ and $\phi_2$ be the two unknown resistances corresponding to the permeation of the capillary walls and the diffusion from the perivascular environment, using the following definitions:
\begin{itemize}
    \item $f_c^{tpz} = I$\, ,
    \item $c_v^{tpz} - c_t^{p,tpz} = \Delta \xi_1$\, ,
    \item $c_t^{p,tpz} - c_t^{tpz} = \Delta \xi_2$\, ,
    \item $c_v^{tpz} - c_t^{tpz} = \Delta \xi_1 + \Delta \xi_2 = \Delta \xi$\, ,
\end{itemize}
they result $\displaystyle \phi_1 = \frac{V}{SP}$ and $\displaystyle \phi_2 = \frac{VL}{SD_t^{tpz}}$\, .
As a consequence, {the drug flux $I$} reads:
\begin{equation}
    I = \frac{\Delta \xi}{\phi_1 + \phi_2} = f_c^{tpz} = \frac{c_v^{tpz}-c_t^{tpz}}{\frac{V}{S}(\frac{1}{P} + \frac{L}{D_t^{tpz}})}
    =\frac{c_v^{tpz}-c_t^{tpz}}{\frac{V}{S}(\frac{D_t^{tpz} + P\cdot L}{PD_t^{tpz}})} =
    \frac{S}{V}\frac{(c_v^{tpz}-c_t^{tpz})P\cdot D_t^{tpz}}{D_t^{tpz} + P\cdot L}\, ,
\end{equation}
where we introduce the constant $K^{tpz} = \frac{S}{V}\frac{P\cdot D_t^{tpz}}{D_t^{tpz} + P\cdot L}$.

Note that the intermediate variable \( c_t^{p,tpz} \) is eliminated algebraically by assuming a quasi-steady transport regime. Under this assumption, the flux of TPZ through the capillary wall and the tissue interface is the same (i.e., the same current \( I \) flows through both resistive elements). This leads to the system:
\[
I = \frac{c_v^{tpz} - c_t^{p,tpz}}{\phi_1} = \frac{c_t^{p,tpz} - c_t^{tpz}}{\phi_2}\, ,
\]
from which one obtains the explicit expression:
\[
c_t^{p,tpz} = \frac{\phi_2 c_v^{tpz} + \phi_1 c_t^{tpz}}{\phi_1 + \phi_2}.
\]
Substituting this relation into the flux expression produces a closed-form relationship between vascular and tissue concentrations, with an effective total resistance \( \phi_1 + \phi_2 \). 

In addition, we consider the coefficient $m$ defined in Eq.~\eqref{metcttpz} as:
\begin{equation*}
    m^{tpz}(c_t^{tpz},c_t^{ox}) =\ \left(\frac{K}{K+c_t^{ox}}\right) \left(k_{met}c_t^{tpz} + \frac{V_{max}{c_t^{tpz}}}{K_{m} + c_t^{tpz}}\right)\, ,
\end{equation*}
and we combine it with the dynamics of the survival fraction described by \eqref{dlogSF}, so that nonlinear dynamics are integrated while conserving a linear functional form.

The reciprocal influence of oxygen and tirapazamine concentrations is also considered:
\begin{equation*}
d_t c_t^{ox} + \phi(SF) m^{ox}(c_t^{ox}) =  P^{ox} \frac{S}{V} (c_v^{ox} - c_t^{p,ox}) + \frac{D_t^{ox}}{L} \frac{S}{V} (c_t^{p,ox}-c_t^{ox})\, ,
\end{equation*}
$c_t^{ox}$ and $c_v^{ox}$ corresponding to tissue and vascular blood oxygen concentrations, while $c_t^{p,ox}$ represents its perivascular concentration and $m^{ox}(c_t^{ox})$ is defined in \eqref{eq:MichaelisMentenConc}.
Exploiting, as before, the electrical analogy between mass and charge transport, we define 
$K^{ox}=\frac{S}{V}\frac{P^{ox} D^{ox}}{D^{ox} + P^{ox} L}$.
In addition, we account for cell viability after drug exposure, including the effect of the surviving fraction. This would account for the impact of cell death, induced by drug concentration, on oxygen consumption rates by integrating pharmacodynamic responses.

As a result, the surrogate model for pharmacodynamics in the tissue governing the time evolution of the variables $c_t^{tpz}$, $c_t^{ox}$, and $SF$ is:
\begin{equation}
\begin{cases}
  d_t c_t^{tpz} + K^{tpz} c_t^{tpz} + \phi_0\,  SF m^{tpz}(c_t^{tpz},c_t^{ox})= K^{tpz} c_v^{tpz}\\
  d_t c_t^{ox} + \phi_0\,  SF\, m^{ox}(c_t^{ox}) + K^{ox} c_t^{ox} = K^{ox} c_v^{ox}\\
 m^{tpz}(c_t^{tpz},c_t^{ox}) {=}\ \left(\frac{K}{K+c_t^{ox}}\right) \left(k_{met}{c^{tpz}_t}+\ \frac{V_{max}{c_t^{tpz}}}{K_{{m}}\ +\ {c_t^{tpz}}}\right) \\
 {m(c^{ox}_t) = V_{max} \frac{c^{ox}_{t}}{c^{ox}_{t} + \alpha^{ox}_t p_{m_{50}}}}\\
\frac{d\ \mathrm{log}\ SF}{dt} = -\alpha \cdot c_t^{tpz} \cdot m^{tpz}(c^{tpz}, c_{t}^{ox}) \\
c_t^{tpz}(0)=0\\
c_t^{ox}(0)=c_0^{ox}
\end{cases}
\label{model_0D}
\end{equation}

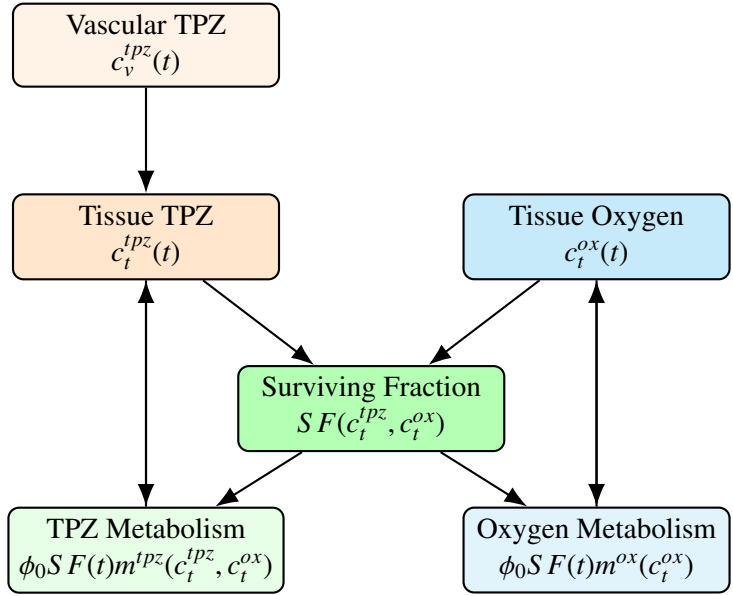
\begin{figure}[ht]
\centering
\begin{tikzpicture}[
  every node/.style={font=\small},
  box/.style={draw, thick, align=center, rounded corners, minimum width=3.5cm, minimum height=1.1cm},
  ellipsebox/.style={draw, ellipse, thick, align=center, minimum width=3.2cm, minimum height=1.0cm},
  arrow/.style={thick, -{Latex[length=3mm]}},
  >=Latex,
  node distance=1.4cm and 2.4cm
]

% Top nodes
\node[box, fill=orange!10] (cv) {Vascular TPZ\\$c_v^{tpz}(t)$};
\node[box, below=of cv, fill=orange!20] (ctpz) {Tissue TPZ\\$c_t^{tpz}(t)$};
\node[box, right=of ctpz, fill=cyan!20] (cox) {Tissue Oxygen\\$c_t^{ox}(t)$};

% Surviving Fraction block
\node[box, below=1.7cm of $(ctpz)!0.5!(cox)$, fill=green!30] (sf) {Surviving Fraction\\$SF(c_t^{tpz}, c_t^{ox})$ };

% Metabolism blocks
\node[box, below=3.0cm of ctpz, fill=green!10] (mtpz) {TPZ Metabolism\\$\phi_0 SF(t)m^{tpz}(c_t^{tpz}, c_t^{ox})$};
\node[box, below=3.0cm of cox, fill=cyan!10] (mox) {Oxygen Metabolism\\$\phi_0 SF(t) m^{ox}(c_t^{ox})$};

% Arrows (Transport and concentration flow)
\draw[arrow] (cv) -- (ctpz);

% Arrows to SF
\draw[arrow] (ctpz) -- (sf);
\draw[arrow] (cox) -- (sf);

% Arrows from SF to metabolisms
\draw[arrow] (sf) -- (mtpz);
\draw[arrow] (sf) -- (mox);

% Arrows to metabolisms
\draw[arrow] (ctpz) -- (mtpz);
\draw[arrow] (cox) -- (mox);

\draw[arrow] (mtpz) -- (ctpz);
\draw[arrow] (mox) -- (cox);

\end{tikzpicture}
\caption{Schematic of the 0D surrogate model. Vascular TPZ drives the tissue concentration \( c_t^{tpz} \), which, together with the tissue oxygen \( c_t^{ox} \), determines the surviving fraction \( SF(t) \). This in turn modulates both TPZ and oxygen metabolism through the nonlinear terms \( m^{tpz} \) and \( m^{ox} \), capturing the feedback structure of the pharmacodynamics model described in equation~\ref{model_0D}.}
\label{fig:0d_model_dependencies}
\end{figure}

Equations \eqref{model_0D},  illustrated in the schematic of Figure \ref{fig:0d_model_dependencies}, define the nonlinear system of ordinary differential equations that constitutes the surrogate model for TPZ pharmacokinetics and pharmacodynamics in tissue complemented with oxygen dynamics. Equation \eqref{model_0D} governs the time evolution of the spatially averaged tissue drug concentration \( c_t^{tpz} \), incorporating both linear extravasation from the vasculature and nonlinear metabolic consumption modulated by oxygen availability and cell viability. Precisely, the dynamics of tissue oxygen concentration \( c_t^{ox} \) is accounted for, where oxygen consumption is also modulated by \( SF(t) \). 

\subsection{One-way interaction between the 0D and the 3D-1D models}
\label{sec:oneway_interaction}

The surrogate model (0D model) presented in Section \ref{sec:model0D} not only provides a standalone framework to investigate pharmacokinetics and pharmacodynamics interactions under the assumption of spatial homogeneity, but also acts as a computationally efficient auxiliary tool within the more complex spatially resolved 3D-1D model. In particular, we employ the 0D model to facilitate and accelerate the evaluation of the nonlinear metabolic consumption terms required by the full multiscale 3D-1D model. This section details the methodology of using a one-way coupling approach to bridge these two modeling scales effectively.

The principal goal of this one-way interaction is to simplify the complex metabolic consumption of TPZ in the 3D-1D framework employing a spatially separable and linear approximation. Specifically, the metabolic consumption is approximated by
\begin{equation}
    \label{eq:linear_metabolism}
    m^{tpz}(\mathbf{x}, t) = SF(t) \cdot r(t) \cdot c^{tpz}(\mathbf{x}, t),
\end{equation}
where $SF(t)$ represents the surviving fraction influenced by tissue concentrations of TPZ and oxygen, and $r(t)$ acts as an effective metabolic rate coefficient integrating both Michaelis–Menten kinetics and oxygen modulation. Thus, the spatial variability of the TPZ concentration, $c^{tpz}(\mathbf{x}, t)$, is preserved.
{The linearization in Eq.~\eqref{eq:linear_metabolism} represents a strategic decomposition of what would otherwise be a computationally prohibitive nonlinear metabolic term encapsulating the complex Michaelis-Menten kinetics with oxygen modulation that must be evaluated pointwise throughout the entire 3D-1D computational domain at each time step. By assuming uniform spatial dependence in the kinetic behavior of metabolic processes, we decomposed this complex function into separable temporal and spatial components: the temporal factors $SF(t)$ and $r(t)$ capture all nonlinear dynamics through pre-computation using the efficient 0D model, while the spatial factor $c^{tpz}(\mathbf{x}, t)$ preserves the essential spatial heterogeneity of drug distribution. This decomposition transforms the original nonlinear problem requiring iterative solving at every grid point into a simple linear operation, dramatically reducing computational cost while maintaining the spatial resolution critical for accurate modeling of drug transport and distribution.}
This linear representation not only circumvents computationally demanding pointwise nonlinear evaluations, significantly enhancing numerical tractability and computational efficiency, but is also directly motivated by implementation constraints. The existing C++ code that underpins the full 3D-1D model illustrated in Figure~\ref{fig:hybrid_model_schematic} is specifically limited to linear reaction terms within the advection-diffusion-reaction solver used to model the dynamics of TPZ. This practical constraint further justifies the adoption of the proposed linear approximation.

\begin{figure}
    \centering
    \begin{minipage}[b]{0.45\textwidth}
        \centering
        \includegraphics[width=\textwidth]{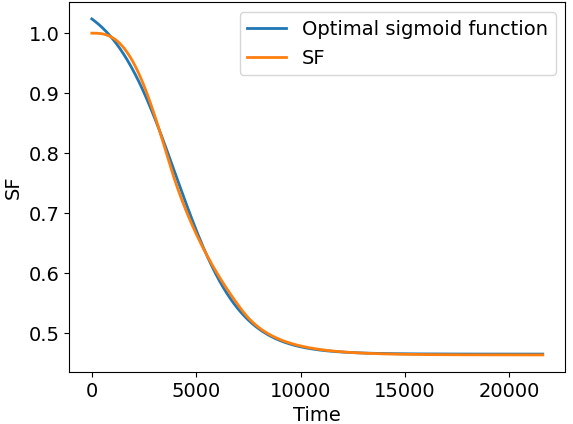}
    \end{minipage}
    \
    \begin{minipage}[b]{0.465\textwidth}
        \centering
        \includegraphics[width=\textwidth]{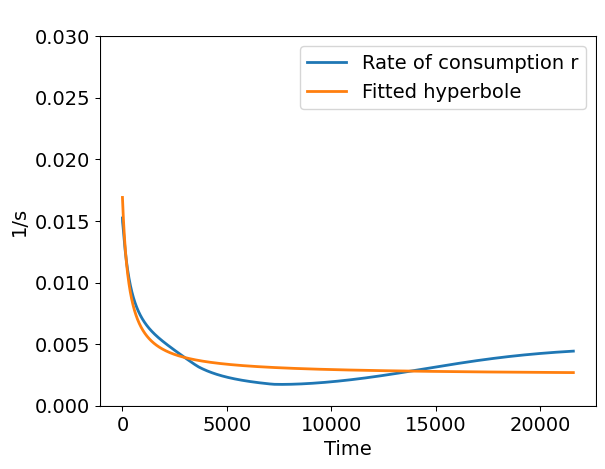}
    \end{minipage}
    \caption{Left panel: Comparison between the surviving fraction ($SF$) computed by the original 0D model and its sigmoid approximation. The remarkable correspondence underscores the adequacy of the sigmoid representation to capture the primary nonlinear transition observed in the 0D simulations. 
    Right panel: Comparison of the computed metabolic rate $r(t)$ from the 0D model and its fitted rational approximation. The figure highlights both the successful capture of the general declining trend and areas where discrepancies remain, potentially indicating more complex underlying dynamics.}
    \label{fig:sigmoid_hyperbole}
\end{figure}

To achieve this efficient representation, we first approximate the surviving fraction $SF(t)$ obtained from the 0D simulations using a sigmoid function:
\begin{equation}
    \label{eq:sigmoid_SF}
    SF(t) = \mathrm{X} - \frac{\mathrm{Y}}{1 + \exp(-\mathrm{Z}(t - D))}.
\end{equation}
The sigmoid approximation is characterized by parameters that represent distinct physiological behaviors: $\mathrm{X}$ and $\mathrm{X}-\mathrm{Y}$ correspond to the upper and lower asymptotes, $\mathrm{Z}$ controls the steepness of the transition, and $D$ marks the midpoint time of the response. The optimized values of these parameters, derived from the regression against the data from the 0D model, are summarized in Table~\ref{tab:sigmoid_SF}. {The parameter optimization was performed in Python using \texttt{scipy.optimize.curve\_fit}, which implements a nonlinear least-squares solver based on the Levenberg--Marquardt algorithm.}
\begin{table}
\centering
\caption{Optimized sigmoid parameters for the surviving fraction $SF(t)$.}
\begin{tabular}{cccc}
\toprule
Parameter & Value & Unit & Description \\
\midrule
$\mathrm{X}$ & 1.0658 & -- & Upper asymptote \\
$\mathrm{Y}$ & 0.6013 & -- & Sigmoid amplitude \\
$\mathrm{Z}$ & 0.00065 & s$^{-1}$ & Slope \\
$D$ & 4002.48 & s & Midpoint time \\
\bottomrule
\end{tabular}
\label{tab:sigmoid_SF}
\end{table}
The high accuracy of this sigmoid fit is quantitatively validated by a coefficient of determination $R^2 = 0.9999$, indicating a satisfactory level of agreement with {data obtained from the 0D model}. However, the Kolmogorov-Smirnov test suggests slight deviations in the residual distribution, reflecting minor systematic discrepancies. This result is further illustrated in Figure~\ref{fig:sigmoid_hyperbole} (left panel), where the fitted sigmoid curve closely follows the original computed $SF(t)$ values of the 0D model.

Following the approximation of $SF(t)$, we compute the effective metabolic coefficient $r(t)$ by dividing the metabolic rate $m^{tpz}(t)$, as output by the 0D model, by the product of $SF(t)$ and the spatially averaged TPZ concentration $c^{tpz}(t)$. The resulting time series for $r(t)$ is suitably approximated using a rational function:
\begin{equation}
    \label{eq:hyperbolic_r}
    r(t) = \frac{\mathrm{A}}{t + \mathrm{B}} + \mathrm{C}\, .
\end{equation}
{This functional form effectively captures the sharp early decay and nonzero asymptotic level observed in the 0D-model data while minimizing residual bias; it did not systematically over- or under-predict data across time. Its hyperbolic structure naturally captures the expected saturation kinetics of metabolic processes, where high initial rates occur when substrate availability is abundant, followed by a monotonic decay to lower steady-state levels as the system approaches metabolic equilibrium. This form provides excellent numerical stability compared to exponential alternatives while maintaining analytical tractability for integration into the 3D-1D framework. Additionally, the three-parameter structure offers interpretable coefficients that correspond to distinct physiological phases: the scaling parameter $\mathrm{A}$ governs the magnitude of initial metabolic activity, the time shift $\mathrm{B}$ controls the transition timing, and the asymptotic offset $\mathrm{C}$ represents the long-term basal metabolic rate. This parsimony in parameters, combined with the function's ability to approximate multi-exponential decay processes typical of complex metabolic networks, makes it an efficient and physiologically meaningful choice for capturing the essential dynamics of TPZ metabolism.}
The fitted parameters, obtained through regression analysis, are listed in Table~\ref{tab:metabolism_fit}.
\begin{table}
\centering
\caption{Fitted parameters for the effective metabolic rate coefficient $r(t)$.}
\begin{tabular}{cccc}
\toprule
Parameter & Value & Unit & Description \\
\midrule
$\mathrm{A}$ & 4.7865 & -- & Scaling numerator \\
$\mathrm{B}$ & 331.6163 & s & Time shift \\
$\mathrm{C}$ & 0.00247 & s$^{-1}$ & Long-term offset \\
\bottomrule
\end{tabular}
\label{tab:metabolism_fit}
\end{table}
The quality of this hyperbolic fit, though moderate with an $R^2 = 0.7000$, demonstrates an adequate ability to capture the main trends and stabilize metabolic decay, even if some residual temporal variability remains unexplained.
{We note that the rational fit of the effective metabolic rate coefficient yielded a coefficient of determination $R^2 \approx 0.70$, which is moderate. However, this level of accuracy is acceptable in the present context because the fit is used as an intermediate surrogate function rather than a direct predictive endpoint. Its main role is to preserve the monotonic oxygen dependence and approximate magnitude of the metabolic rate, both of which are reliably captured.} Figure~\ref{fig:sigmoid_hyperbole} (right panel) depicts both the original metabolic rate data and the fitted rational function, visually illustrating the strengths and limitations of the approximation.

In summary, substituting equations~\eqref{eq:sigmoid_SF} and~\eqref{eq:hyperbolic_r} into~\eqref{eq:linear_metabolism} produces a highly efficient and analytically explicit expression for the metabolic source term in the 3D-1D model. This combined approach significantly improves computational performance, thus facilitating extensive parametric investigations, optimization studies, and uncertainty quantification without sacrificing spatial resolution.

\begin{figure}[ht]
\centering
\begin{tikzpicture}[
  every node/.style={font=\small},
  box/.style={draw, rounded corners, thick, align=center, minimum width=4.0cm, minimum height=1.0cm},
  arrow/.style={thick, -{Latex[length=3mm]}},
  doublearrow/.style={thick, <->, >=Latex},
  >=Latex,
  node distance=1.0cm and 1.0cm
]

% Top layer: Vascular network
\node[box, fill=blue!25] (vascular) {Vascular Network (1D)\\[1pt] Blood Flow \& Hematocrit};

% Middle left: Oxygen transport
\node[box, below left=of vascular, fill=cyan!35] (oxygen) {Oxygen Transport (3D-1D)\\[1pt] Diffusion, Exchange};

% Middle right: Drug delivery
\node[box, below right=of vascular, fill=orange!35] (drug) {TPZ Delivery (3D-1D)\\[1pt] Transport, Activation};

% Bottom center: 0D
\node[box, below=of $(drug)!0.5!(oxygen)$, fill=green!15] (0D) {0D model\\[1pt] $SF(t)$};

% Lower left: Oxygen in tissue
\node[box, below=of oxygen, fill=cyan!20] (ctox) {Tissue Oxygen\\$c_t^{ox}(t,\mathbf{x})$};

% Lower right: Drug in tissue
\node[box, below=of drug, fill=orange!20] (ctpz) {Tissue TPZ\\$c_t^{tpz}(t,\mathbf{x})$};

% Bottom center: SF
\node[box, below=of $(ctox)!0.5!(ctpz)$, fill=green!15] (sf) {Cell Survival (0D)\\[1pt] $SF(t)$};

% Bottom center: r
\node[box, below=of sf, fill=green!15] (r) {TPZ consumption (0D)\\[1pt] $SF(c^{tpz}, c_{t}^{ox})$};

% Bottom center: mox
\node[box, below=of ctox, fill=cyan!10] (mox) {Oxygen metabolization\\[1pt] $\phi_0 SF(t) m^{ox}(c^{ox})$};

% Bottom center: mtpx
\node[box, below=of ctpz, fill=orange!10] (mtpz) {TPZ metabolization\\$\phi_0 SF(t) r(t) c_{t}^{tpz}$};

% Arrows between boxes
\draw[arrow] (vascular) -- (oxygen);
\draw[arrow] (vascular) -- (drug);
\draw[arrow] (oxygen) -- (ctox);
\draw[arrow] (drug) -- (ctpz);
\draw[arrow] (0D) -- (sf);
\draw[arrow] (ctox) -- (mox);
\draw[arrow] (ctpz) -- (mtpz);
\draw[arrow] (sf) -- (mox);
\draw[arrow] (sf) -- (mtpz);
\draw[arrow] (sf) -- (r);
\draw[arrow] (mox) -- (ctox);
\draw[arrow] (mtpz) -- (ctpz);
\draw[arrow] (r) -- (mtpz);

\end{tikzpicture}
\caption{Schematic representation of the multiscale 3D-1D-0D pharmacokinetics model architecture after implementing the one-way interaction with the 0D model. Tissue-level TPZ and oxygen concentrations are computed via the 3D transport equations, while the corresponding metabolic source terms are no longer evaluated through nested nonlinear functions. Instead, TPZ metabolism is modeled as a linear expression modulated by two surrogate functions derived offline from the 0D model: the surviving fraction \( SF(t) \) and the effective metabolic coefficient \( r(t) \). Oxygen metabolism is similarly represented via an exogenous effective function \( m_{\mathrm{ox}}^{\text{eff}}(t) \). This reformulation reduces computational complexity while preserving the essential physiological feedback.}
\label{fig:hybrid_model_schematic_0D}
\end{figure}
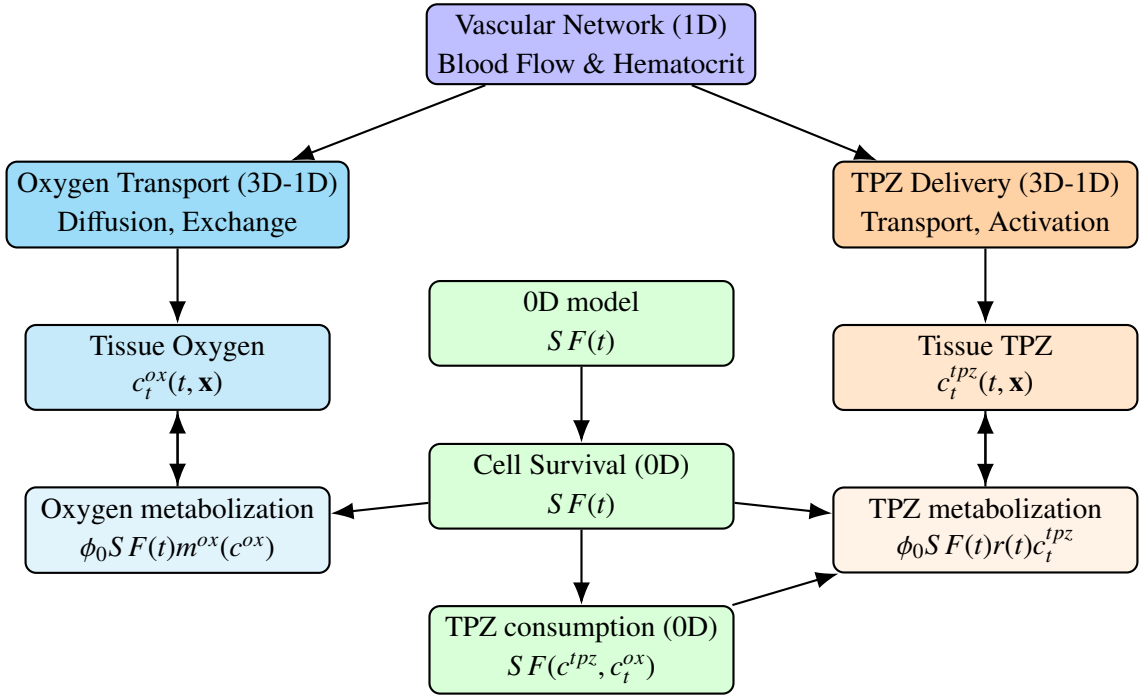

The conceptual structure of the updated 3D-1D model, incorporating the surrogate functions \( SF(t) \), \( r(t) \), and {an exogenous effective function \( m_{\mathrm{ox}}^{\text{eff}}(t) \)}, is summarized in Figure~\ref{fig:hybrid_model_schematic_0D}. We refer to this reformulated system as the 3D-1D-0D model, to emphasize the multiscale architecture that combines spatially resolved 3D and 1D transport dynamics with surrogate functions derived from a reduced 0D pharmacokinetics/pharmacodynamics model. This framework will serve as the basis for the numerical experiments of Section \ref{sec:results}. 

{Note that, in this paper, the classical \textit{full-memory} formulation is adopted (see Eq.~\eqref{defSF}). As such, the history of the system is fully taken into account through the integral operators. An alternative approach consists of using \textit{short-memory} kernels, in which the influence of distant past states progressively vanishes or is truncated after a finite horizon. Such models, recently investigated in the context of viscoelasticity and continuum mechanics~\cite{ActaMech2020, IJSS2020}, may provide a more realistic description of biological processes characterized by fading memory and could also reduce computational cost. Although not considered in the present study, the extension of our framework to short-memory operators represents an interesting direction for future work.} {In particular, fractional-order models could be employed to account for memory effects in tissue transport and drug–cell interactions, potentially capturing delayed cellular responses or heterogeneous diffusion behaviors not represented in classical integer-order models. Incorporating such effects may improve model fidelity in settings where drug uptake or clearance exhibits non-exponential kinetics. Exploring these generalizations lies beyond the scope of the present study but represents an exciting avenue for future research.}

\subsection{Mathematical analysis of the 0D and the 3D-1D pharmacokinetics models}

{The analysis of the existence of solutions for the proposed coupled multiscale model presents significant theoretical challenges due to the nonlinear couplings between flow dynamics, transport phenomena, and pharmacodynamic responses. Although a complete existence proof for the full coupled system is beyond the scope of this study, we provide here a mathematical framework that establishes the foundation for such an analysis.}

{Our multiscale model consists of four main components whose existence properties can be analyzed separately:
\begin{itemize}
\item ($\mathcal{F}\&\mathcal{H}$): The microvascular flow model is based on well-established Poiseuille flow in vessels and Darcy flow in tissue. For systems of this type, the results of existence and uniqueness are available when the hydraulic conductivities and geometric parameters meet the standard regularity conditions~\cite{quarteroni2016cardiovascular}. The hematocrit transport Eq.~\eqref{eq:1D_H} represents a hyperbolic conservation law with source terms, for which existence follows from standard theory~\cite{bressan2000hyperbolic}.
\item $\mathcal{O}$: The coupled 3D-1D oxygen transport system in \eqref{eq:Ofull} consists of advection-diffusion-reaction equations with Michaelis-Menten kinetics. The existence of weak solutions for such systems is well-established in the literature~\cite{evans2010partial}. 
\item $\mathcal{TPZ}$: The main mathematical challenge lies in the strong nonlinear coupling through the surviving fraction $SF(t)$ defined in Eq.~\eqref{defSF}. This integral equation creates a global-in-time dependence that requires careful analysis.
\item The 0D model (Section~\ref{sec:model0D}) serves as a regularization of the full system and provides smooth approximations that can be rigorously analyzed. For this simplified system, existence follows directly from the standard ODE theory~\cite{hartman2002ordinary}.
\end{itemize}
Lastly, while formal existence proofs remain as future theoretical work, our extensive numerical experiments provide strong computational evidence of well-posedness.
As such, while rigorous existence analysis represents an important theoretical challenge that merits dedicated mathematical investigation, the current study focuses on computational implementation and sensitivity analysis. The numerical evidence strongly supports well-posed behavior, and the mathematical framework outlined above provides the foundation for future theoretical work.}

\section{Numerical Discretization Techniques}

Given the complexity of the coupled mathematical models that describe the microvascular environment and drug transport, analytical solutions are not available. Therefore, numerical simulations are essential to apply these models to realistic scenarios. This section outlines the numerical techniques employed for the mixed-dimensional 3D-1D models and the 0D surrogate model presented in Section~\ref{mvasMod}.

\subsection{3D-1D Model Discretization}

The core of the spatially distributed model is a \textit{mixed-dimensional 3D-1D} framework, which describes the tissue environment as a three-dimensional (3D) domain and embeds the microvascular network within it as a collection of one-dimensional (1D) channels or a metric graph, see for example~\cite{Cattaneo2014b, Possenti2019, Possenti2021}. The \textit{finite element method (FEM)} is used to discretize the governing partial differential equations (PDEs) for 3D-1D problems, including blood flow, oxygen transfer, and drug transport. This method is based on the variational formulation and the partitioning of the domain into finite elements. We refer the interested reader to specific papers on the formulation and discretization of these equations, for example~\cite{D'Angelo2008, D'Angelo2012, Koppl2018953, Laurino2019, Kuchta2021558, Heltai20232425}.
A key advantage of the mixed-dimensional formulation is that the discretizations of the equations defined in the tissue and vascular networks are entirely independent. {Meshes ultimately lack shared topological characteristics (such as vessel segments located at boundaries or centroids of 3D elements), and geometrically consistent features are not necessary in the design of numerical schemes.}  
The tissue is discretized using a uniform tetrahedral mesh. Piecewise continuous polynomial finite elements are used for quantities such as oxygen and drug concentrations, whereas mixed finite elements are used for interstitial fluid flow. The resolution of the mesh is determined through a mesh sensitivity analysis. For example, a typical domain size $500\,\mu\mathrm{m} \times 500\,\mu\mathrm{m} \times 500\,\mu\mathrm{m}$ is discretized with 15 nodes per side.
The 1D branches of the vascular network are discretized as separate subdomains, approximated by straightly segmented pieces, typically divided into five equispaced elements per branch. Continuous piecewise polynomial finite elements are also used for variables such as blood flow and drug transport in the vascular system.
To solve the coupled problem, a linearization strategy is employed to handle nonlinearities, using either a fixed-point iteration or the Newton-Raphson method.  
{The resulting linear systems arising from $\mathcal{F}\&\mathcal{H}$ as well as from $\mathcal{O}$ and $\mathcal{TPZ}$ at each iteration are solved using a direct SuperLU solver.}
All 3D-1D simulations are performed using an \textit{in-house C++ code} built on the open-source \texttt{GetFEM++} library, which enables discretization and coupling of operators across multiple dimensions and supports non-matching grids between embedded and embedding domains.

{From the numerical approximation standpoint, the reduced regularity of the solution induced by line sources leads to suboptimal convergence, with order $O(h^{1/2-\varepsilon})$ in the $H^1$ seminorm on quasi-uniform meshes, being $h$ the mesh characteristic size and $\varepsilon > 0$ an arbitrarily small constant. Nevertheless, numerical analysis of related models with inclusions~\cite{Koppl2018953, Gjerde2019} shows that optimal finite element rates ($O(h)$ in $H^1$) are recovered when the mesh is aligned with the embedded network or suitably refined near the singular sources. In particular, Köppl et al.~\cite{Koppl2018953} provide comprehensive analysis for second-order elliptic problems with inclusions, showing that embedded lower-dimensional structures create singular behavior that requires appropriate mesh treatment for optimal convergence; while, Gjerde et al.\cite{Gjerde2019} show for elliptic equations with line sources that the solution exhibits reduced regularity near the line source and that optimal convergence can be restored by appropriate treatment of singular terms. However, the satisfaction of such conditions for optimal approximation significantly increases the computational cost of the method, making its application to realistic problems (as the one addressed here) computationally intractable.}

\subsection{0D Model Discretization and Integration}

A spatially averaged \textit{surrogate model (0D)} is used to simulate non-linear pharmacokinetics and pharmacodynamics responses (e.g., drug metabolism rate and cell survival fraction $SF$) that depend on quantities at the tissue level, such as oxygen concentration. This model provides input-output curves (e.g., $SF(t)$ and $r(t)$ ) that are parameterized and integrated into the larger 3D-1D FEM framework.
Internally, the 0D model is formulated as a system of \textit{ ordinary differential equations (ODE)}, which account for the time evolution of drug and oxygen concentrations and the nonlinear dependence of $SF$ on those quantities. These ODEs reflect Michaelis-Menten-type kinetics and exponential decay laws for drug action.
The ODE system is numerically solved using MATLAB's \texttt{ODE Suite}, specifically the \texttt{ode45} solver, which implements a Runge-Kutta (4,5) method with adaptive time-stepping. This choice ensures computational efficiency and robustness for the stiff, nonlinear behavior typical of pharmacokinetics and pharmacodynamics models. {Under smoothness assumptions, this scheme achieves a global order of accuracy equal to four, with local error control of order five.}
The results of these simulations, that is, the functions $SF(t)$ and $r(t)$ {defined in Eq.s~\eqref{eq:sigmoid_SF} and~\eqref{eq:hyperbolic_r}, respectively}, are then used to inform the source terms in the 3D-1D FEM framework. This multiscale strategy allows the 3D model to incorporate complex, nonlinear cellular-level dynamics without solving the full set of coupled equations at every spatial node and time step.

\section{Sensitivity Analysis of the Multiscale Model}\label{sec:SA}

As detailed in the preceding section, we have developed a multiscale model to simulate the pharmacokinetics and pharmacodynamics of hypoxia-activated drugs in the vascular microenvironment. 
This model integrates various interconnected components, including blood flow, oxygen transport, and drug delivery, each governed by a set of physical, physiological, and geometrical parameters. 
Given the complexity of this multiscale model and the inherent uncertainty in the precise values of these parameters within the heterogeneous tumor microenvironment, it becomes crucial to assess the robustness of the model predictions and to identify the most influential parameters affecting treatment outcomes, such as drug concentration, oxygen levels, and possibly cell survival.
Therefore, this section presents the methodology adopted for global sensitivity analysis, a technique essential to systematically exploring how variations in the input parameters of our multiscale model impact its outputs across their physiological and pathological ranges. 
By identifying the most influential factors, this analysis aims to provide valuable insight into the design and optimization of hypoxia-targeted cancer therapies and guide future experimental investigations. 
This section provides a concise overview of the methodology used to perform a global sensitivity analysis that assesses the influence of variations in input across the full spectrum of potential input values.

\subsection{Variance based methods}

{Variance-based sensitivity analysis methods, such as Sobol’s indices, are widely regarded as the gold standard for global sensitivity assessment, as they decompose the output variance into contributions from individual parameters and their interactions~\cite{Saltelli2008}. These approaches provide rigorous quantitative information but require a very large number of model evaluations, which is prohibitive for computationally expensive models such as the current 3D-1D framework. For this reason, variance-based methods are not employed in this work; instead, we focus on the Morris screening approach, which offers an efficient compromise between interpretability and computational feasibility, while still allowing us to capture the main effects and potential interactions among parameters.}

\subsection{Screening methods: Elementary effect}

The elementary effects (EE) method is simple but effective in screening for a few important input factors over the many that can be contained in a model. 
The fundamental idea behind the method was proposed by Morris in 1991 with the definition of the concept of elementary effects~\cite{morris1991factorial}.
The EE method determines whether an input factor is \textit{negligible}, \textit{linear and additive}, \textit{ non-linear or interacts with some other factor}~\cite{saltelli2004global}. This test corresponds to an average of derivatives over the space of inputs. 
{Let $X\in \mathbb{R}^k$ represent the vector of input factors for a model and $Y \in \mathbb{R}$ its scalar outcome. By assuming the $k$ components of $X$ as independent factors, each normalized to the interval $[0,1]$, the input space is defined by the $k$-dimensional unit hypercube $\mathcal{X}_k:=[0,1]^k$. For the construction of elementary effect trajectories, the hypercube $\mathcal{X}_k$ is discretized on $p\in\mathbb{N}$ levels, therefore inducing a partition $\mathcal{X}_k^p$ of the parameter domain $[0,1]^k$ from which admissible candidate points are drawn.
For a given input $X$, the elementary effect corresponding to the $i$-th component $X_i$, with $i=1,\ldots,k$ is defined as:
\begin{equation*}
EE_{i}(X)={\frac{[Y(X+e_i\Delta)-Y(X)]}{\Delta}}
\end{equation*}
being $\Delta$ a value in $\displaystyle \left\{\frac{1}{p-1},\frac{2}{p-1},\cdots,\frac{p-2}{p-1}\right\}$ and $e_i$ the unitary vector for the $i$-th direction. 
It is important to mention that for each $i=1,\dots,k$, $EE_i(X)$ captures local characteristics since it relies on $X$. To derive a broader sensitivity measure, one should determine the statistics of their distribution, denoted as $F_i\sim EE_i$, by randomly selecting various $X$ from $\mathcal{X}_k^p$.
Therefore, after sampling $r$ points $X^{(1)},\dots,X^{(r)}$, we obtain the following sensitivity metrics for the $i$-th component of the input factors:
%Note that for each $X\in\mathcal{X}_k^p$ it is required that $X+e_i\Delta\in\mathcal{X}_k^p$ for $i=1,\cdots,k$. The distribution $F_i$ of the $i$-th input factor is obtained by collecting the elementary effects computed along $r$ one-at-a-time trajectories within the discretized hypercube $\mathcal{X}_k^p$, where each consecutive point differs from the previous one in only one factor. A measure of the sensitivity of $Y$ with respect to $X_i$ is then given by the mean of these $r$ realizations:
\begin{equation}
\mu_{i}=\frac{1}{r}\sum_{j=1}^rEE_i(X^{(j)})\, .
\end{equation}
We observe that the index $\mu_{i}$, representing the sensitivity of the $i$-th parameter globally, can be deceptive when faced with non-monotonic relationships, due to the potential offsetting of positive and negative contributions.
For this reason, we introduce \(\mu^*_i\), the mean of the absolute values of the elementary effects, which more reliably quantifies the overall importance of each parameter irrespective of the direction of influence:
\begin{equation}
    \mu^*_{i}=\frac{1}{r}\sum_{j=1}^r|EE_i(X^{(j)})|
\end{equation}
The values of $\mu^*_i$ are obtained as the mean of the absolute values of the elementary effects computed along the sampling points within $\mathcal{X}_k^p$.
The $EE_i$ standard deviation $\sigma_i$ is used to estimate the quality of the effect, namely, whether it results from nonlinear effects or due to mutual interactions with other factors:
\begin{equation*}
\sigma_{i}=\sqrt{\frac{1}{r-1}\sum_{j=1}^r\left(EE_i(X^{(j)})-{\mu}_{i}\right)^{2}}\,.
\end{equation*}
Calculating each $EE_i$ requires two sample points, which leads to $2rk$ evaluations of the model in order to calculate the measures $\mu^*_i$ and {$\sigma_i$ for $i=1,\ldots,k$}.} 

{To reduce computational cost, Morris proposed a unique sampling technique that allows the use of only $r(k+1)$ samples \cite{morris1991factorial}. This involves $r$ trajectories, each consisting of $k+1$ points in the unit hypercube $\mathcal{X}_k^p$ that differ from their adjacent point in a single component, thus yielding $k$ elementary effects for each trajectory. Each trajectory can be interpreted as a matrix ${\bf B}^*$ of size $(k+1)\times k$, constructed in such a way that each input variable is varied one at a time while all others remain fixed. Specifically, for each $j=1,\dots,k$, there is a pair of rows in which only the $j$-th component differs. 
However, since ${\bf B}^*$ is built in a fixed and deterministic way, the resulting variations would not reflect the randomness required for a proper global sensitivity analysis. To restore randomness, two additional components are used:
(i) a diagonal matrix $\mathbf{D}^*$ of size $k \times k$ with diagonal entries randomly chosen from $\{-1, 1\}$, used to assign a random sign (positive or negative) to each step in the trajectory;
(ii) a permutation matrix $\mathbf{P}^*$ of size $k \times k$, used to randomly change the order in which the input variables are perturbed.
The diagonal matrix $\mathbf{D}^*$ ensures that a variable can be increased or decreased by $\Delta$ in the trajectory. The permutation matrix $\mathbf{P}^*$ ensures that the order of these changes is not fixed, adding further variability to the sampling process.
By evaluating the model at each point of the resulting perturbed trajectory, one can calculate an elementary effect for each input variable as the difference between two consecutive model evaluations, divided by $\Delta$. These differences provide a local sensitivity measure of the model output with respect to each input. For a more detailed discussion, we refer the reader to \cite{Saltelli2008}.}

\section{Results of the Hypoxia-Activated Drug Model and Sensitivity Analysis}\label{sec:results}

This section presents the results obtained by coupling the proposed multiscale mathematical model with global sensitivity analysis techniques to evaluate the influence of key parameters on model outcomes. The primary goal is to identify the most significant physiological, microvascular, and drug-related parameters that affect therapeutic efficacy, with a particular focus on the cell survival fraction ($SF$).
To this end, a series of numerical simulations was carried out, incorporating a time-dependent tirapazamine (TPZ) injection profile, described in detail in Section~\ref{sec:TPZprofile}. The model tracks the TPZ concentration in the tissue ($c_t^{tpz}$) and the corresponding $SF$ as indicators of drug distribution and treatment effectiveness.
Given the large number of parameters and their wide variability, direct sensitivity analysis on the full model would be computationally prohibitive. Therefore, a two-stage approach was adopted. A preliminary screening was conducted using a surrogate formulation (0D model) to identify the most influential parameters. The sensitivity analysis was then refined on the detailed 3D-1D-0D model, focusing only on this reduced subset of parameters.

Although physiological bounds for most parameters were available in the literature, the data for others were incomplete or lacking. In such cases, a variability range of $\pm25\%$ around the physiological baseline was assumed. Table~\ref{tab:boundsSA} summarizes the ranges of parameters and their respective physiological bounds.
\begin{table}
\centering 
\footnotesize
\begin{tabular}{|c|c|c|c|c|c|}
\hline
i & $X_i^{min}$ & $X_i^{max}$ & Physiological Bounds  \\ \hline \hline
$ c_{v0}^{tpz}$~\cite{senan1997phase}  & $1.78\times 10^{-2} \, \, mol/m^3$ & $4.73\times 10^{-2} \, \, mol/m^3$ & $ 200 - 330\, \,  mg/m^2 $ \\ \hline
$ D_t^{tpz}$~\cite{hicks1998extravascular} & $1.80\times 10^{-11} \, \, m^2/s$ &$1.25\times 10^{-10} \, \, m^2/s$ & $ 0.18 - 1.25 \times 10^{-6} cm^2/s$ \\ \hline
$ P_{tpz}$~\cite{Vitullo2023} & $ 3.75 \times 10^{-5} \, \, m/s$&$6.25\times 10^{-5} \, \, m/s$ & $ 5\times 10^{-5}\, \, m/s \pm 25\% $ \\ \hline
$ k_{met}$~\cite{hicks1998extravascular} & $ 5.00\times 10^{-3} \, \, s $ & $3.33\times 10^{-2} \, \, s $ & $ 0.3 - 2\, \, min $ \\ \hline
$ V_{max}^{tpz}$~\cite{PKPD} & $ 1.07\times 10^{-4} \, \, mol/(m^3\, s) $ & $ 1.78\times 10^{-4} \, \, mol/(m^3\, s) $ & $ 1.42 \times 10^{-4}\, \, mol/(m^3\, s) \pm 25\% $ \\ \hline
$ K_m^{tpz}$~\cite{PKPD} & $ 2.63\times 10^{-3} \, \, mol/m^3 $ & $ 4.38\times 10^{-3} \, \, mol/m^3 $ & $ 3.5 \times 10^{-3} \, \, mol/m^3 \pm 25\% $ \\ \hline
$ K $~\cite{kyle1999measurement} & $ 2.60\times 10^{-3} \, \, mol/m^3 $ & $ 1.30\times 10^{-2} \, \, mol/m^3 $ & $ 2 -10\, \, mmHg $ \\ \hline
$ \alpha $~\cite{PKPD} & $ 1.75\times 10  \, \, (mol/m^3)^{-2} $&$ 2.91\times 10 \, \, (mol/m^3)^{-2} $ & $ 23.3\, \,  (mol/m^3)^{-2} \pm 25\% $ \\ \hline
$ \phi_0 $~\cite{PKPD} & $ 3.88\times 10^{-1} $ & $ 6.46\times 10^{-1} $ & $ 0.517 \pm 25\% $ \\ \hline
$ c_{v0}^{ox} $ & $ 3.90\times 10^{-2} \, \,  mol/m^3 $ & $ 1.30\times 10^{-1} \, \,  mol/m^3 $ & $ 30 - 100mmHg $ \\ \hline
$ V_{max}^{ox} $~\cite{Vitullo2023,CZ} & $ 1.30\times 10^{-3} \, \, mol/m^3/s  $ & $ 1.04\times 10^{-2} \, \, mol/m^3/s $ & $ 1 - 8\, \,  mmHg/s $  \\ \hline
$ K_m^{ox}$~\cite{Vitullo2023} & $ 6.50\times 10^{-4} \, \, mol/m^3/s $ & $ 1.30\times 10^{-3} \, \, mol/m^3/s $ & $ 0.5 - 1 mmHg $ \\ \hline
$ P_{ox} $~\cite{Vitullo2023} & $ 3.50\times 10^{-5} \, \, m/s $ & $ 3.00\times 10^{-4} \, \, m/s $ & $ - $ \\ \hline
$ D_{ox} $~\cite{Possenti2021} & $ 1.81\times 10^{-9} \, \, m^2/s $ & $ 3.01\times 10^{-9} \, \, m^2/s $ & $ 2.41\times 10^-9\, \, m^2/s \pm 25\%  $ \\ \hline
\end{tabular}
\caption{Parameter ranges for the sensitivity analysis along with the relative physiological bounds.}
\label{tab:boundsSA}
\end{table}

\subsection{{Synthetic Vascular Network Generation}}

{The vascular networks adopted in this study are generated according to a biomimetic and iterative procedure originally proposed in~\cite{Possenti2019}, designed to satisfy key physiological and morphological criteria. The construction relies on Voronoi tessellations to model capillary topologies that ensure space-filling properties while preserving biological realism. Each network is obtained by generating planar Voronoi tessellations from a random distribution of $8$ seed points on a square domain with $500\,\mu$ m side, resulting in networks with a sufficiently high aspect ratio (\( L / R \approx 4 \)) for each branch to justify the one-dimensional approximation of blood flow.
To assign physiologically consistent vessel radii, an iterative refinement algorithm is used, enforcing Murray's law in bifurcations and anastomoses. The initial radii are set uniformly to $4\, \mu$m and are iteratively updated based on flow-driven topological connectivity. For each bifurcation, a random split ratio $\frac{R_{\text{in},0}}{R_{\text{out},1}}$ is assigned, and the radii of the daughter vessels are determined by the following rules:
(i) bifurcation, $R_{\text{in},0}^3 = R_{\text{out},1}^3 + R_{\text{out},2}^3$, 
(ii) anastomosis, $R_{\text{in},1}^3 + R_{\text{in},2}^3 = R_{\text{out},0}^3$.
This procedure is repeated until convergence. To maintain physiological significance, we keep only those configurations where radii fall within the range $[2, 6]$ $\mu$m and the average radius is equal to $\overline{R} = 4 \, \mu$m with a tolerance of $\pm$ $5\%$. }

\begin{figure}
    \centering
    \includegraphics[width=0.5\linewidth]{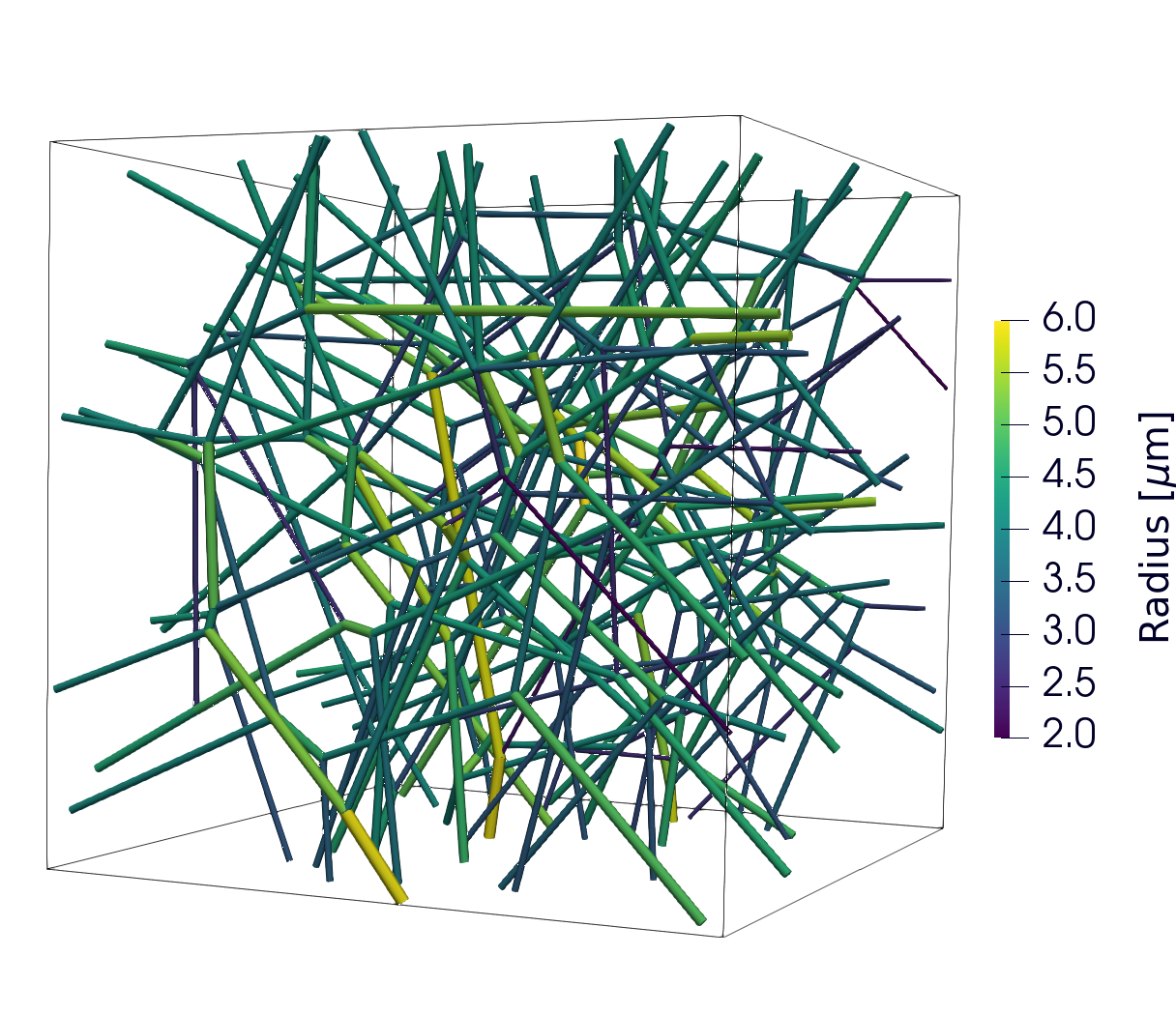}
    \caption{Example of constrained Voronoi-based synthetic vascular network and visualization of radii distribution ($\mu m$). This is the vascular network that has been used consistently in all the numerical tests.}
    \label{fig:voronoi_net}
\end{figure}

{We construct a collection of approximately \( 10^4 \) valid configurations.
3D vascular networks are then obtained by stacking $18$ admissible planar networks along the vertical axis, with each layer occupying a slab of thickness \( T_z = 27.8 \, \mu\text{m} \). To avoid unrealistic planarity, the nodes are vertically perturbed within the range \(\pm T_z/2\), ensuring a physiologically plausible spatial distribution. The resulting 3D network achieves a surface-to-volume ratio of approximately \(\approx 7000 \, \text{m}^{-1} \)~\cite{Baxter1990a}, which corresponds to the target physiological value.
This structured approach allows controlled variability through a small number of parameters, while ensuring compliance with essential physiological and topological constraints. All the results of this work are obtained using a single synthetic vascular network with an associated radius distribution reported in Figure~\ref{fig:voronoi_net}.}

\subsection{Computational Setup and TPZ Injection Profile}
\label{sec:TPZprofile}

The three-dimensional tissue domain was modeled as a cubic sample of size $500\,\mu\text{m} \times 500\,\mu\text{m} \times 500\,\mu\text{m}$ and discretized using a uniform tetrahedral mesh with 15 nodes per edge. Linear finite elements were used for spatial discretization, resulting in approximately $4358$ degrees of freedom (DOF) in the tissue domain.
The discretization of the one-dimensional vascular network was performed based on the specific topology of the embedded vessel architecture. Specifically, we have used $5021$ DOFs to discretize a problem with $181$ vessels. 
A mesh sensitivity analysis was performed to ensure that spatial resolution was sufficient to capture the relevant features of oxygen and drug transport, without introducing unnecessary computational overhead (see also~\cite{Vitullo2023}).
The computational cost of simulating the 3D-1D-0D model is highly dependent on the complexity of the vascular network, the mesh density, and the domain size. For example, a single simulation could range from tens of minutes to a few hours on a standard workstation, depending on the number of vessels and the simulation time window. These substantial computational requirements justify the development and use of reduced-order models, such as the 0D surrogate formulation introduced in Section~\ref{sec:model0D}, particularly for large-scale parametric analyses such as global sensitivity studies.

{For the sensitivity analysis based on the Morris method, we set $p = 10$ and $r = 70$. Although it differs from standard selection ($p = 4$, $r = 10$, as noted in \cite{Saltelli2008}), this particular choice was made to meet the requirements of our multiscale model. A finer discretization ($p = 10$) was used to better resolve localized nonlinear effects and was paired with a sufficiently high number of trajectories ($r = 70$) to ensure that the increased resolution translated into meaningful statistical estimates. This choice was informed by previous experience \cite{Vitullo2023}, where we applied a similar methodology to oxygen transport in the vascular microenvironment and found that although ranking stability was achieved around $r\approx 20$, full convergence of the sensitivity indices required $r\approx 50$. For the current study, we opted for $r = 70$ to ensure robustness, even at the expense of an increase in computational cost. Finally, the perturbation step was calculated using the standard formulation $\Delta= p/[2(p-1)]\approx 0.55$, ensuring valid finite-difference approximations within the unit hypercube and consistency with previous implementations.}

The outputs of interest included $c_t^{tpz}$ and $SF$ taken at three specific time points: end of infusion ($t = 7200\,\text{s}$, or $t = 2\,\text{h}$), mid-decay ($t = 10800\,\text{s}$, or $t = 3\,\text{h}$), and end of simulation ($t = 21600\,\text{s}$, or $t = 6\,\text{h}$), as well as their averaged time values $\overline{c}_t^{tpz}$ and $\overline{SF}$ throughout the observation period.
TPZ perfusion is simulated using a time-dependent boundary condition on vascular concentration $c_v$, mimicking a typical chemotherapy treatment. The injection protocol consists of three distinct phases:

\begin{itemize}
    \item[$i)$] Infusion phase ($0\,\text{s} < t \leq 3600\,\text{s} \ | \ 0\,\text{h} < t \leq 1\,\text{h}$): {$c_{v, in}^{tpz}$} increases linearly, simulating the gradual introduction of the drug.
    
    \item[$ii)$] Sustained Infusion Phase ($3600\,\text{s} < t \leq 7200\,\text{s} \ | \ 1\,\text{h} < t \leq 2\,\text{h}$): {$c_{v,in}^{tpz}$} is kept constant, representing steady drug administration.
    
    \item[$iii)$] Post-infusion Decay Phase ($7200\,\text{s} < t \leq 21600\,\text{s}  \ | \ 2 \,\text{h} < t \leq 6\,\text{h}$): {$c_{v, in}^{tpz}$} undergoes exponential decay to model drug clearance, with a characteristic time constant $\tau = 3220\,\text{s}$ ($\simeq 54 \ \text{min}$)~\cite{TPZchildren}.
\end{itemize}
The proposed dynamic injection methodology facilitates an accurate simulation of drug kinetics, enabling a comprehensive evaluation of therapeutic effectiveness by analyzing the progression of $c_t^{tpz}$ and $SF$ across critical time points and throughout the treatment.

\subsection{Sensitivity Analysis for the 0D Model}

To systematically identify the most influential parameters that govern the pharmacokinetics and pharmacodynamics of the 0D model, we apply the Morris method of Elementary Effects (EE), presented in Section \ref{sec:SA}. The results are summarized using two key statistical indices: the mean of the absolute values of the elementary effects, denoted by \(\mu^*\), and their standard deviation, \(\sigma\).  
The parameters with low \(\mu^*\) and low \(\sigma\) can be considered negligible; high \(\mu^*\) and low \(\sigma\) indicate strong, nearly linear effects; while high values of both \(\mu^*\) and \(\sigma\) suggest non-linear or interaction-driven influences. 
A commonly adopted visualization plots each parameter in the $\mu^*$-$\sigma$ plane. This type of plot offers an intuitive representation to discern not only which parameters matter most but also how their effects manifest in the model's behavior.
\begin{figure}
    \centering
    \includegraphics[width=0.5\linewidth]{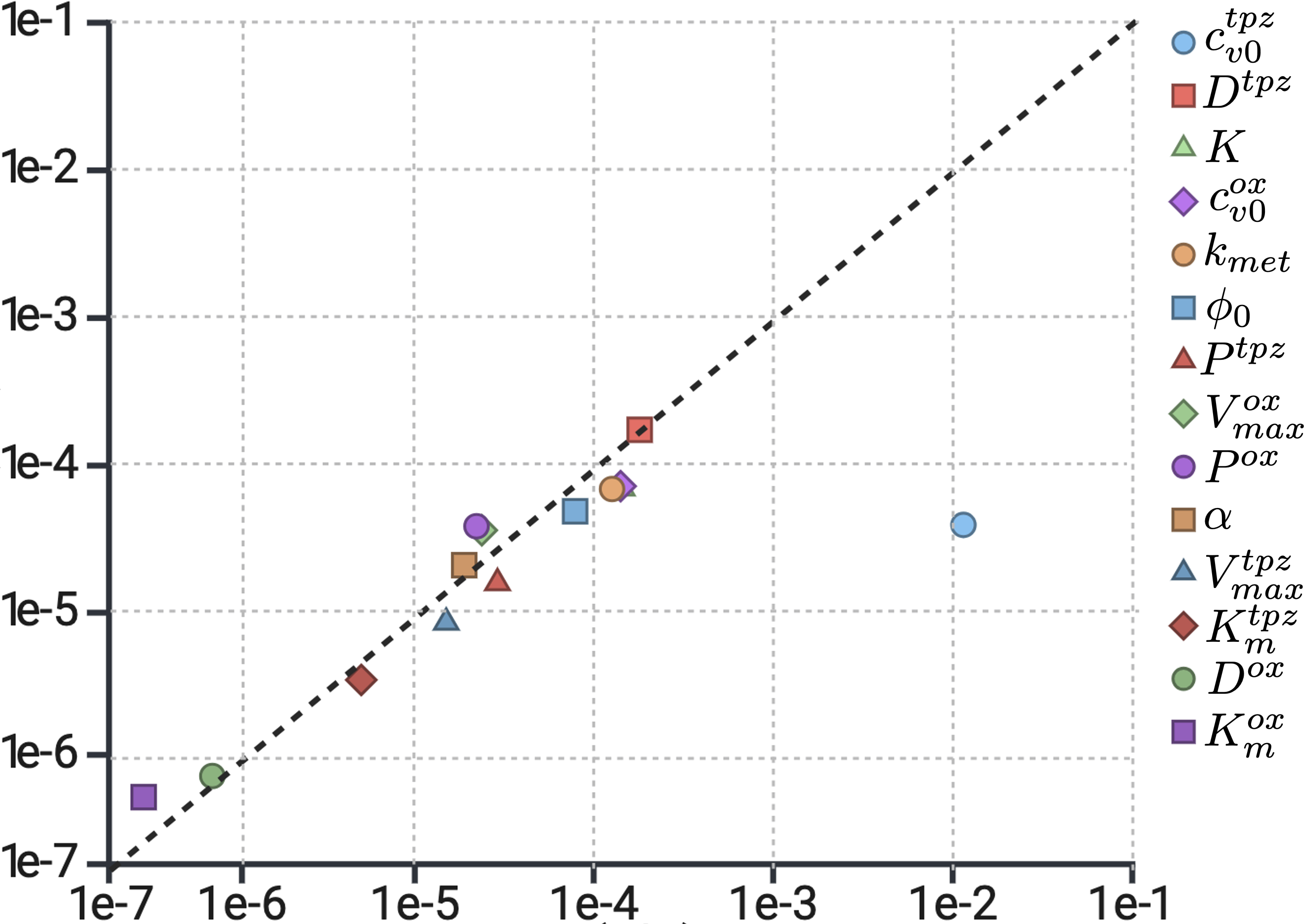}
    \caption{Morris sensitivity indices represented in the $\mu^*$–$\sigma$ plane for the output 
$\overline{c}_t^{tpz}$ obtained with the 0D model. The horizontal axis reports the mean of the absolute elementary effects ($\mu^*$), measuring the overall influence of each parameter, while the vertical axis reports their standard deviation ($\sigma$), which quantifies non-linear effects and interactions. }
    \label{ctaverage}
\end{figure}
In this study, this analysis is used to prioritize the role of physiological and pharmacological parameters in determining tissue-level drug concentration \(c^{tpz}_t\) and the surviving fraction \(SF\), helping guide the subsequent refinement of high-fidelity 3D-1D simulations.
Figure~\ref{ctaverage} shows the $\mu^*_i$-$\sigma_i$ plot for $i =$ $c_{v0}^{tpz}$, $D^{tpz}$, $K$, $c_{v0}^{ox}$, $k_{met}$, $\phi_0$, $P^{tpz}$, $V_{max}^{ox}$, $P^{ox}$, $\alpha$, $V_{max}^{tpz}$, $K_m^{tpz}$, $D^{ox}$, and $K_m^{ox}$ for the quantity of interest $\overline{c}_t^{tpz}$. 
We observe that the value of $\mu^*_{c_{v0}^{tpz}}$ is significantly higher, approximately on the order of $10^{-2}$, compared to the value for other parameters, the latter ranging from $10^{-7}$ to $10^{-4}$. 
This calculation enables a direct comparison of the $\mu^*_i$ index with the quantity of interest considered, suggesting that the vascular concentration of TPZ exerts the most significant influence on the average tissue concentration with an average effect of $\sim 10^{-2}$. 
Moreover, $\sigma_{c_{v0}^{tpz}}$ shows relatively modest values compared to $\mu^*_{c_{v0}^{tpz}}$. $\sigma_i$ assesses the combined effects of factors, including both nonlinear effects and interactions with other factors. Low values of $\sigma_i$ suggest minimal variability between elementary effects, implying that the influence of $c_{v0}^{tpz}$ is largely independent of the values assumed by the other factors.
\begin{figure}
\centering
\begin{adjustbox}{width=\textwidth}
\begin{tikzpicture}
\begin{axis}[
    ybar,
    bar width=7pt,
    ymin=1e-9,
    ymode=log,
    log origin=infty,
    ylabel={$ \mu_i^*(c^{tpz}_t)$},
    symbolic x coords={
        $c_{v0}^{tpz}$, $D^{tpz}$, $K$, $c_{v0}^{ox}$, $k_{met}$, $\phi_0$,
        $P^{tpz}$, $V_{max}^{ox}$, $P^{ox}$, $\alpha$, $V_{max}^{tpz}$,
        $K_m^{tpz}$, $D^{ox}$, $K_m^{ox}$
    },
    xtick=data,
    x tick label style={rotate=45, anchor=east},
    enlarge x limits=0.05,
    x=35,
    bar shift auto,
    legend style={at={(0.5,1.05)}, anchor=south, legend columns=-1},]

% Values at 7200 s
\addplot+ coordinates {($c_{v0}^{tpz}$,2.94e-2) ($D^{tpz}$,3.58e-4) ($K$,2.45e-4) ($c_{v0}^{ox}$,2.48e-4) ($k_{met}$,2.68e-4) ($\phi_0$,1.65e-4) ($P^{tpz}$,5.86e-5) ($V_{max}^{ox}$,3.56e-5) ($P^{ox}$,2.70e-5) ($\alpha$,5.69e-5) ($V_{max}^{tpz}$,1.81e-5) ($K_m^{tpz}$,2.17e-6) ($D^{ox}$,1.05e-6) ($K_m^{ox}$,3.15e-7)};

\addlegendentry{7200 s}

% Values at 10800 s
\addplot+ coordinates {
    ($c_{v0}^{tpz}$,2.06e-2) ($D^{tpz}$,3.34e-4) ($K$,2.64e-4) 
    ($c_{v0}^{ox}$,2.58e-4) ($k_{met}$,2.46e-4) ($\phi_0$,1.38e-4)
    ($P^{tpz}$,5.36e-5) ($V_{max}^{ox}$,4.54e-5) ($P^{ox}$,4.26e-5)
    ($\alpha$,2.93e-5) ($V_{max}^{tpz}$,2.25e-5) ($K_m^{tpz}$,4.53e-6)
    ($D^{ox}$,1.25e-6) ($K_m^{ox}$,4.32e-7)
};
\addlegendentry{10800 s}

% Values at 21600 s
\addplot+ coordinates {
    ($c_{v0}^{tpz}$,3.37e-4) ($D^{tpz}$,7.19e-6) ($K$,5.79e-6) 
    ($c_{v0}^{ox}$,5.42e-6) ($k_{met}$,1.76e-6) ($\phi_0$,4.12e-6)
    ($P^{tpz}$,1.09e-6) ($V_{max}^{ox}$,8.44e-7) ($P^{ox}$,6.48e-7)
    ($\alpha$,1.67e-6) ($V_{max}^{tpz}$,2.02e-6) ($K_m^{tpz}$,2.44e-6)
    ($D^{ox}$,2.27e-8) ($K_m^{ox}$,6.55e-9)
};
\addlegendentry{21600 s}

\end{axis}
\end{tikzpicture}
\end{adjustbox}
\caption{Bar chart of $\mu_i^*(c_t^{tpz})$ taken at $7200$, $10800$, and $21600$ seconds, obtained with the 0D-model.}
\label{fig:mu_ct_0D}
\end{figure}
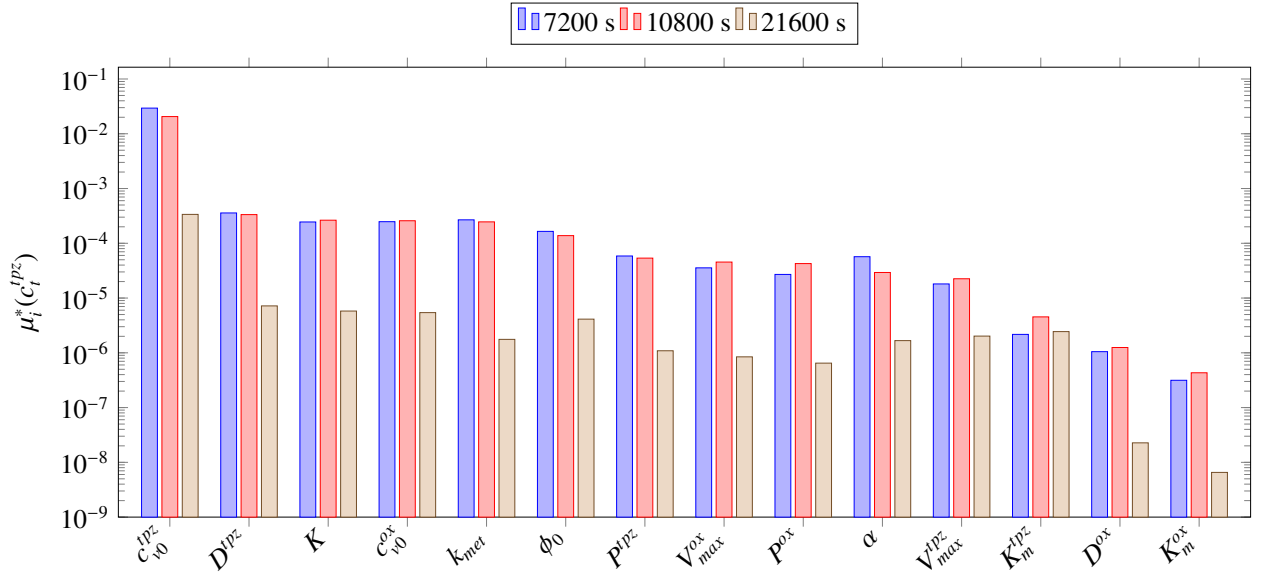

{Note that the quantity $\overline{c}_t^{tpz}$ is used here as a representative of the biomedical problem, being a measure of the amount of chemotherapeutic delivered to the tissue. In fact, the metric $\overline{c}_t^{tpz}$ serves as a preliminary evaluation measure. It is suitable to identify the parameters that most influence bulk tissue exposure, and hence the average pharmacodynamics response, in our scenarios. However, when the clinical or biological question requires the prediction of the fate of small, highly hypoxic niches, a spatially resolved analysis is necessary, and the averaged metric should be replaced by point-wise or near-vessel concentration diagnostics~\cite{Wilson2007, Hicks2010, Rockwell2013}. However, it may not be sufficiently descriptive as a reference quantity for sensitivity analysis.}

\begin{figure}[H]
\centering
\includegraphics[width=0.33\linewidth]{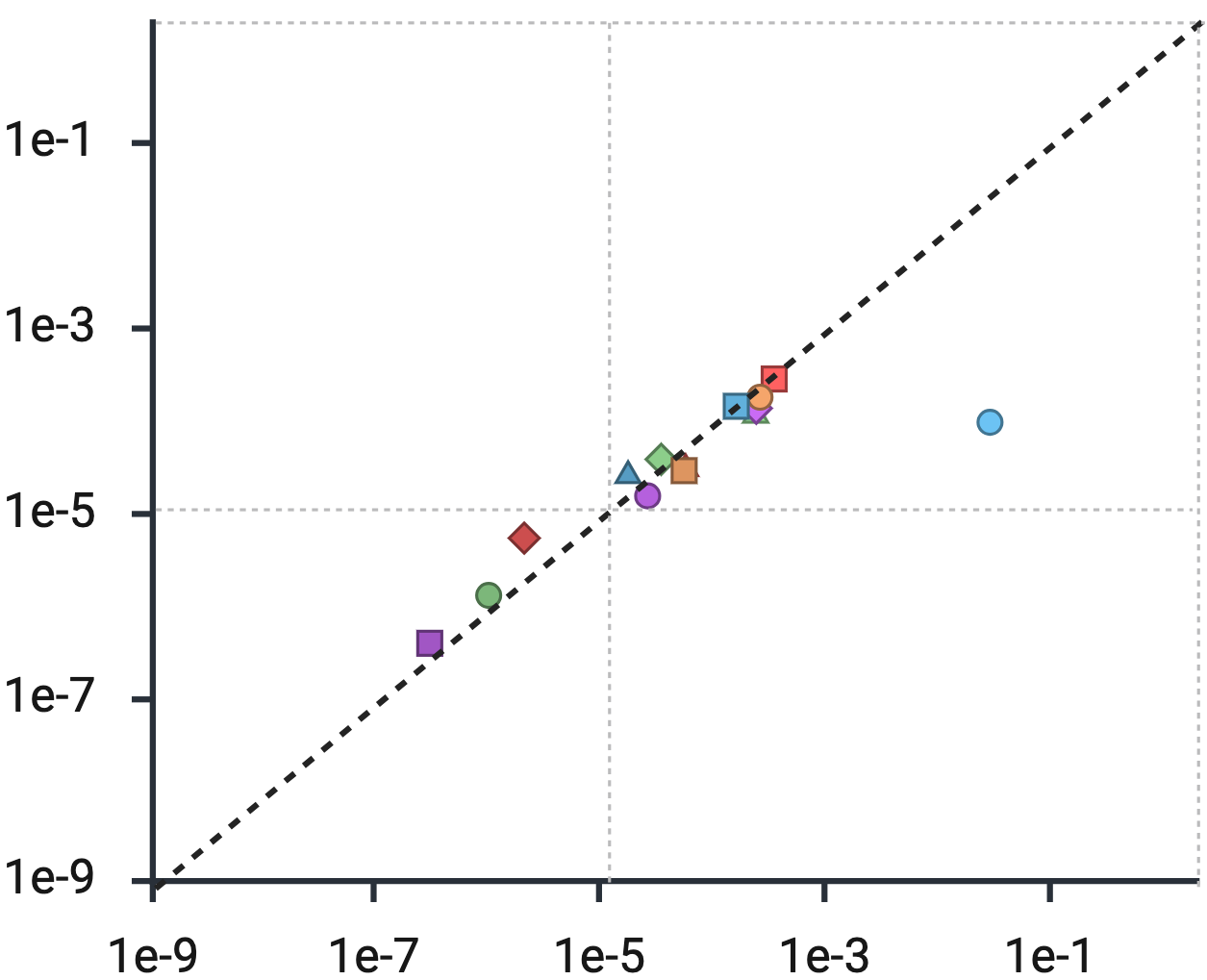}
\hspace{-0.2cm}
\includegraphics[width=0.33\linewidth]{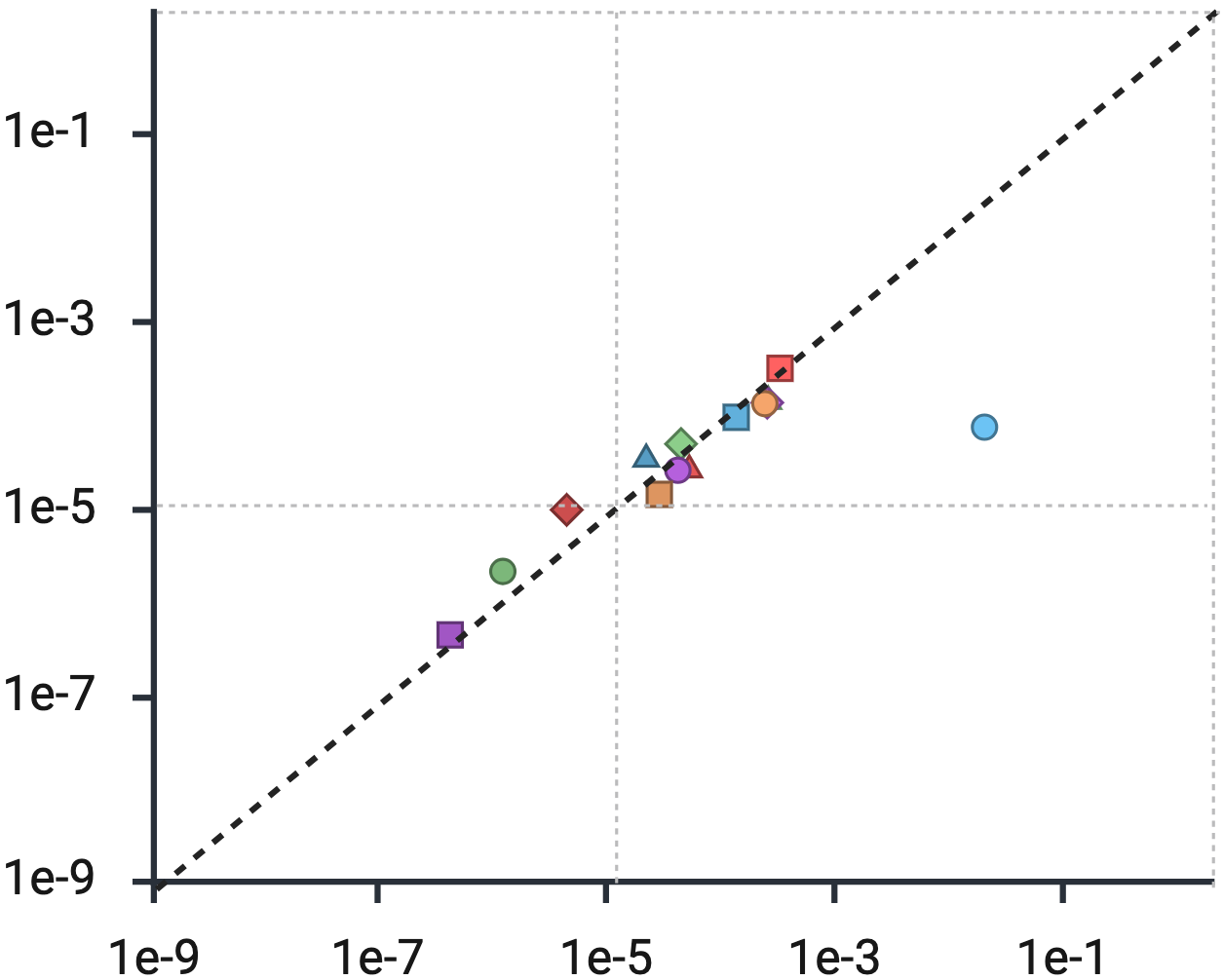}
\hspace{-0.2cm}
\includegraphics[width=0.33\linewidth]{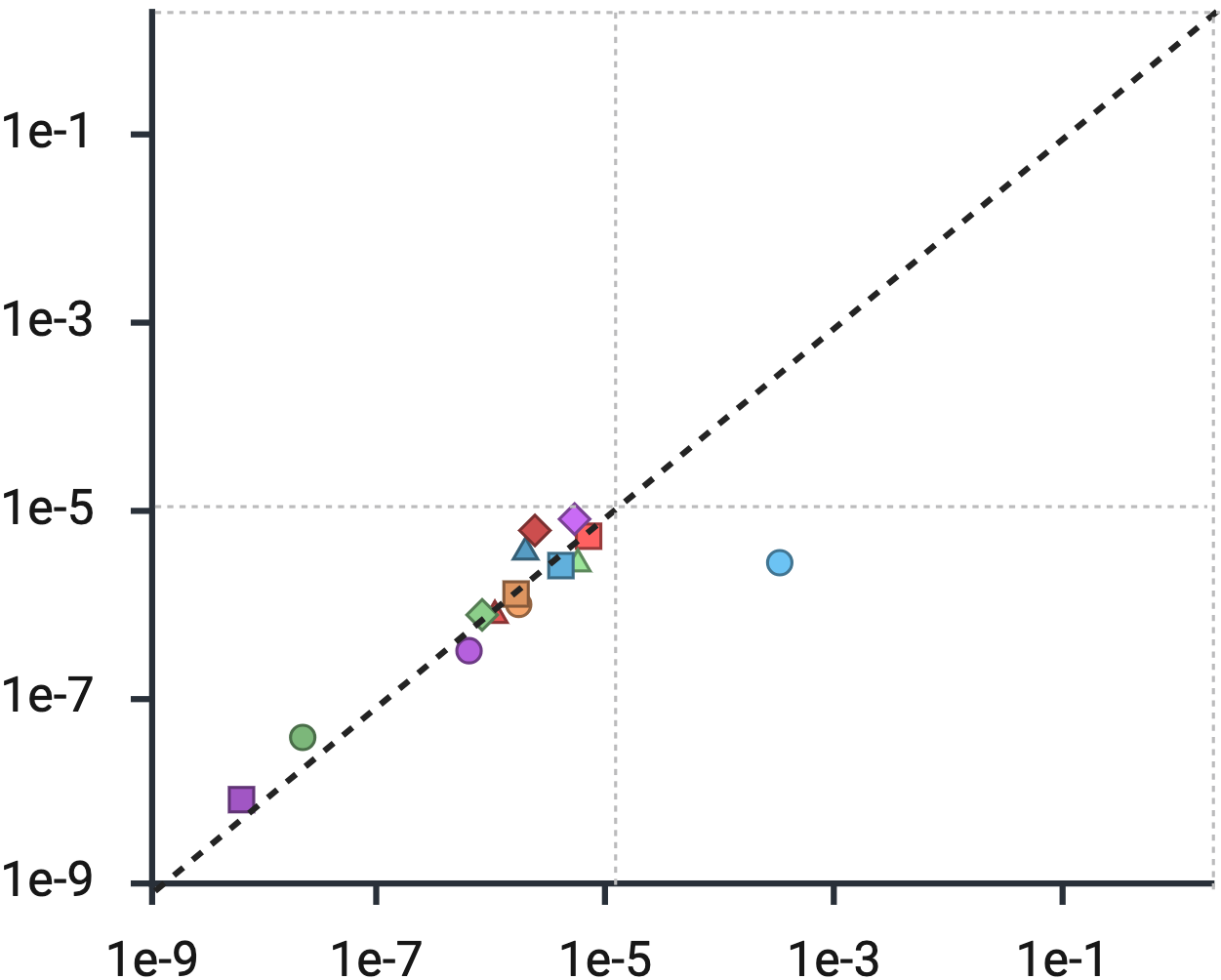} \\
\includegraphics[width=0.9\linewidth]{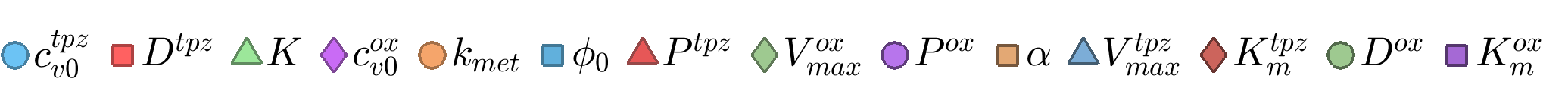}
\caption{Morris indices in the $\mu^*$-$\sigma$ plane 
(the horizontal axis shows $\mu^*$ 
and the vertical axis shows $\sigma$) 
taken at $7200$, $10800$, and $21600$ seconds, relatively to $c_t^{tpz}$, obtained with the 0D-model.}
\label{fig:sigma_ct_0D}
\end{figure}
Due to the exponential decay for $t>7200\, s$, $c_t^{tpz}$ assumes very low values for most of the time of numerical experiments compared to its value at the peak of the infusion phase ($3600\, s < t \leq 7200\, s)$. In this light, the concentration of TPZ in tissue at $t =7200\, s$, $10800\, s$, and $21600\, s$ may provide a more complete picture of the sensitivity of the system in terms of $\mu^*_i$ and $\sigma_i$. As we can observe in Figure~\ref{fig:mu_ct_0D}, interestingly, all parameters have a stronger influence on $c_t^{tpz}(7200\, s)$ and $c_t^{tpz}(10800\, s)$ since the exponential decay of $c_v^{tpz}$ prescribed at the inflow boundary does not influence the results. Moreover, as confirmation, the values of $\mu^*_i$ for $c_t^{tpz}(21600\, s)$ are, indeed, significantly smaller than the other three quantities.
Figure~\ref{fig:sigma_ct_0D} collects Morris indices for $c_t^{tpz}(7200\, s)$, $c_t^{tpz}(10800\, s)$ and $c_t^{tpz}(21600\, s)$. Interestingly, except for $c_{v0}^{tpz}$, all indices lie about on the bisector of the $\mu^*_i$-$\sigma_i$ plane, thus suggesting nonlinear effects of the input parameters on the output of the model and highlighting possible interactions between them.

\begin{figure}[H]
    \centering
    \includegraphics[width=0.5\linewidth]{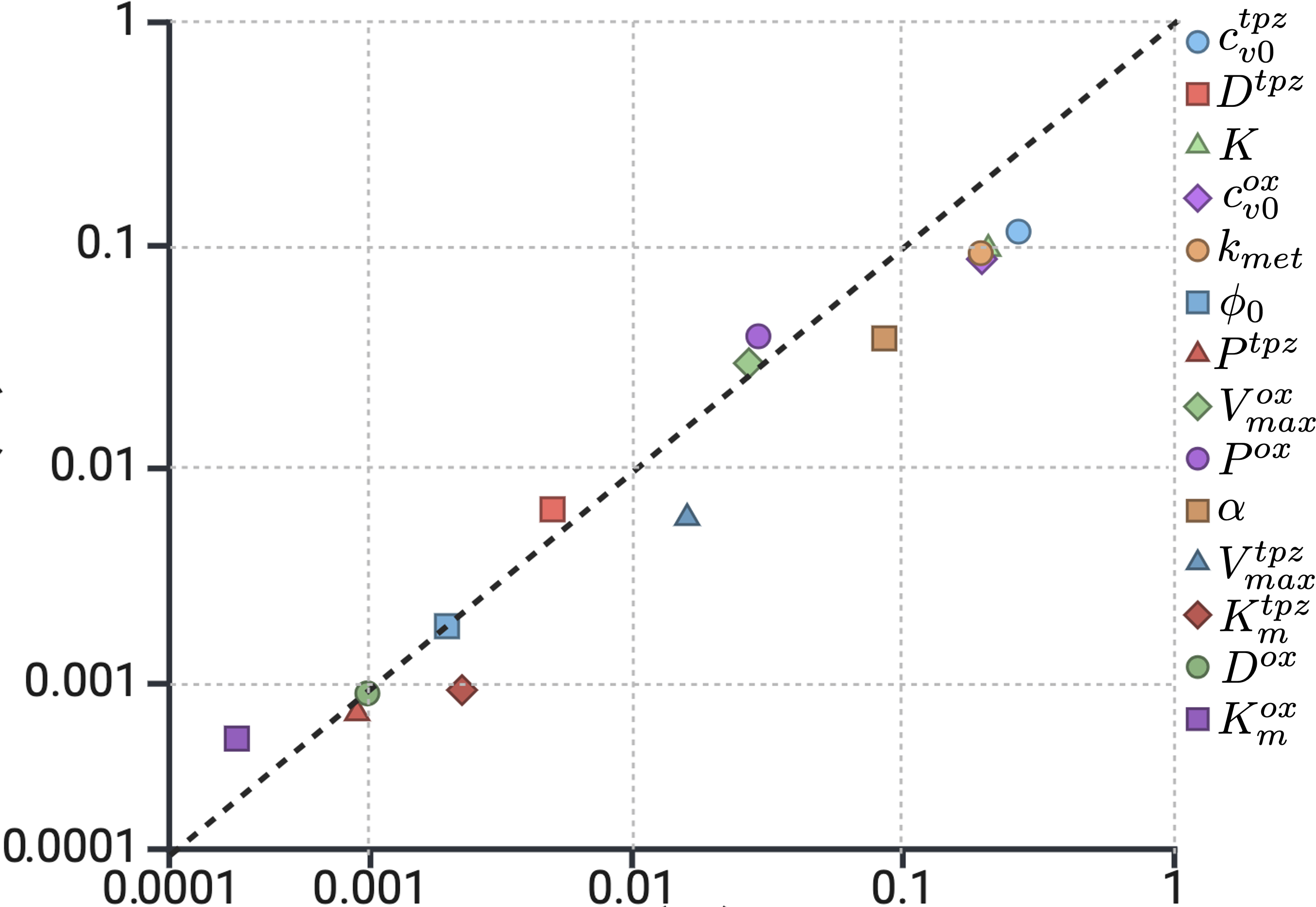}
    \caption{Morris indices in the $\mu^*$-$\sigma$ plane (the horizontal axis shows $\mu^*$ and the vertical axis shows $\sigma$) relative to $\overline{SF}$ obtained with the 0D-model.}
    \label{fig:mu_SF_av}
\end{figure}

The surviving fraction $SF$ is used as a measure of therapeutic efficacy, assessing cell death caused by drug toxicity. 
Similarly to $c_t^{tpz}$, we see that regarding $\overline{SF}$ the most influential factor is $c_{v0}^{tpz}$ (see Figure~\ref{fig:mu_SF_av}). Interestingly, the metabolic parameters also greatly affect the model outputs, particularly those related to oxygen concentration, such as $K$ and $c_{v0}^{ox}$.  
The values of $\sigma_i$ indicate possible interactions between the parameters and the resulting nonlinear effects. 
In general, $c_{v0}^{tpz}$, $k_{met}$, $K$, $\alpha$, $c_{v0}^{ox}$, $V_{max}^{ox}$, and $P_{ox}$ are identified as the most critical parameters for $\overline{SF}$.
\begin{figure}[H]
\centering
\includegraphics[width=0.33\linewidth]{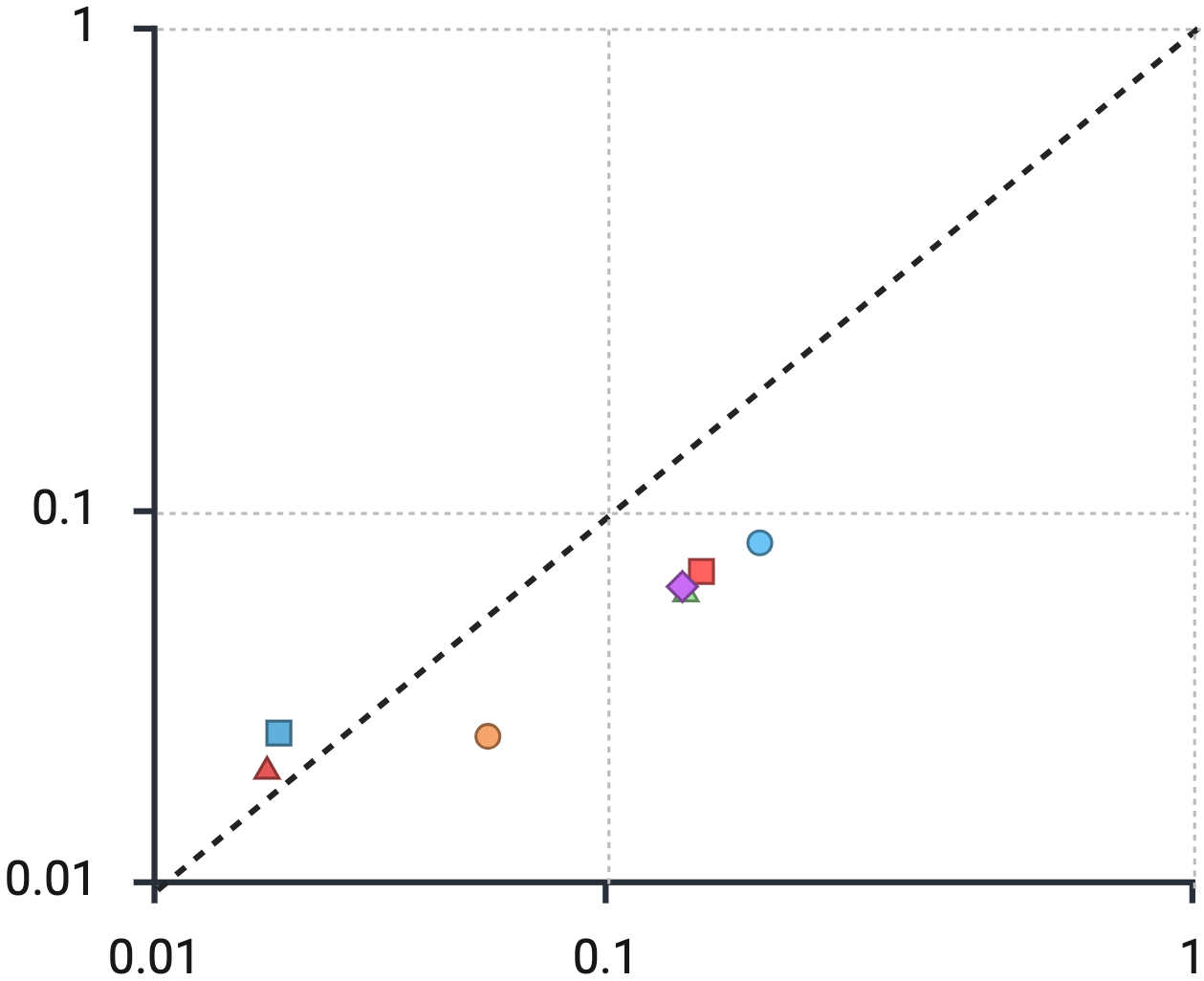}
\hspace{-0.2cm}
\includegraphics[width=0.33\linewidth]{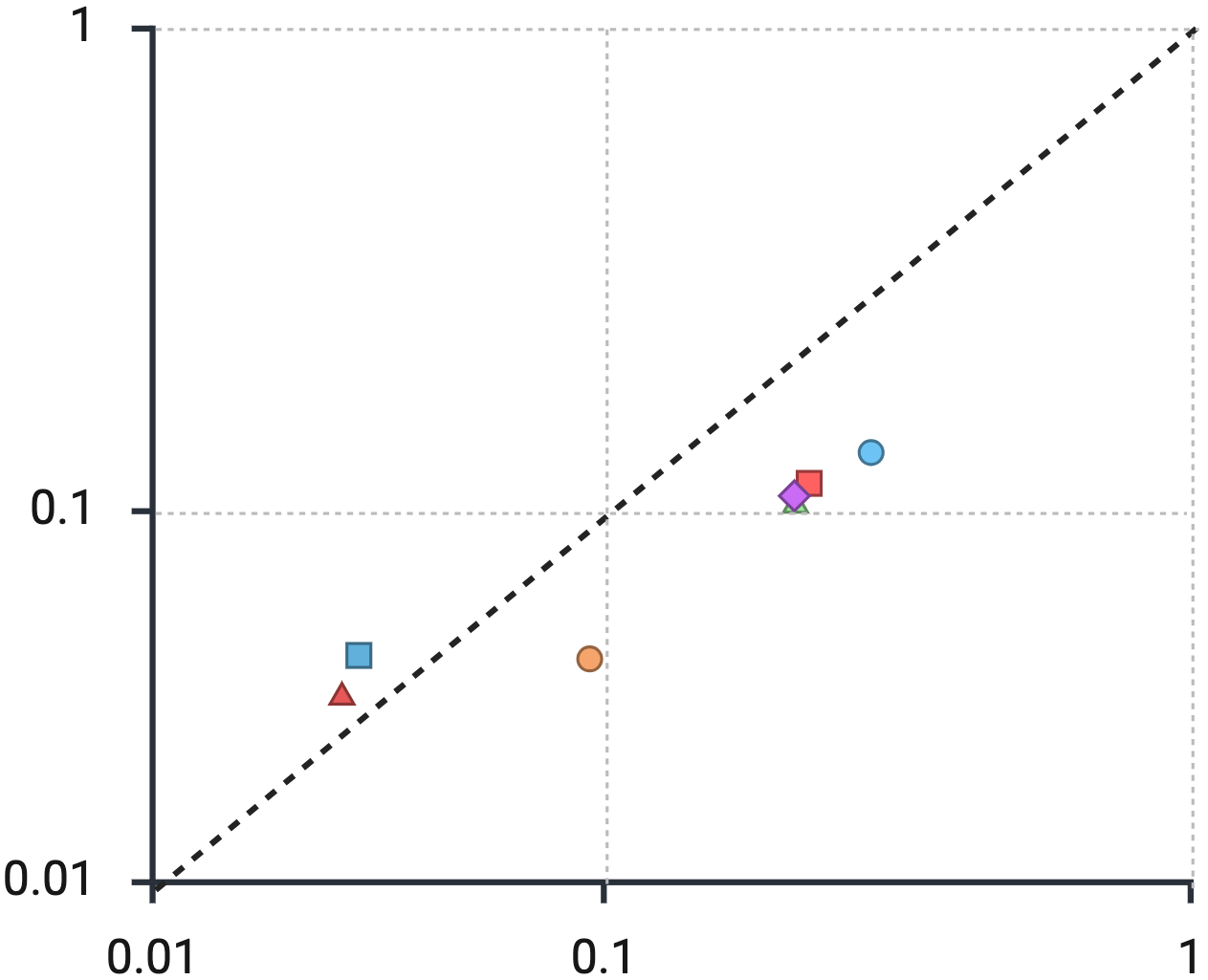}
\hspace{-0.2cm}
\includegraphics[width=0.33\linewidth]{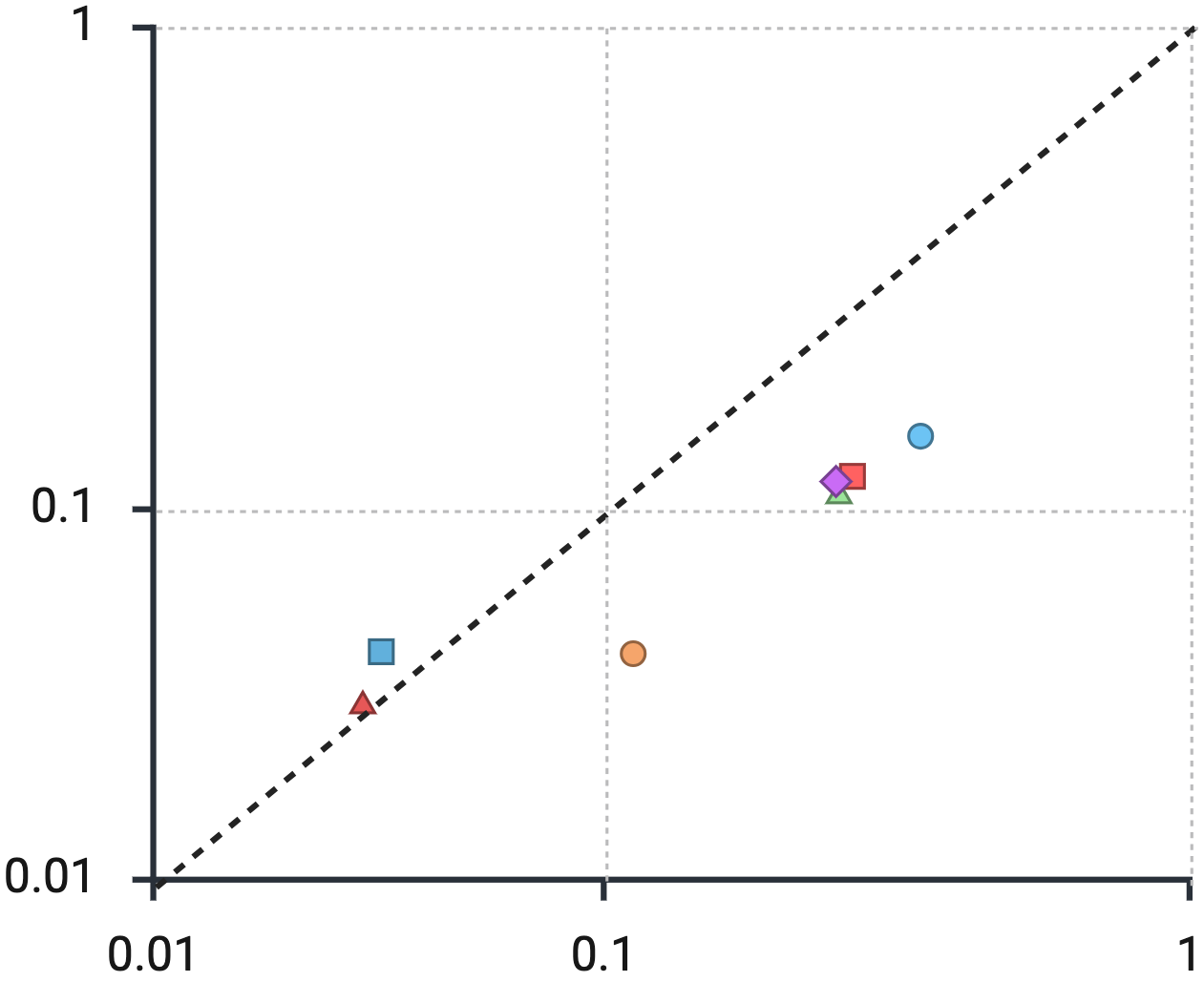} \\
\includegraphics[width=0.5\linewidth]{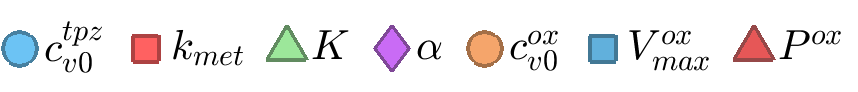}
\caption{Morris indices in the $\mu^*$-$\sigma$ plane (the horizontal axis shows $\mu^*$ and the vertical axis shows $\sigma$) taken at $7200$, $10800$, and $21600$ seconds, relatively to $SF$, obtained with the 0D-model.}
\label{fig:SF_0D}
\end{figure}
Figure~\ref{fig:SF_0D} reports Morris indices for these critical parameters taken at $7200\, s$, $10800\, s$, and $21600\, s$. 
As expected, exposure to TPZ significantly influences model output.
Specifically, focusing on $c_{v0}^{tpz}$,  $\mu^*_i(SF(21600\, s))$ is higher than those taken at previous times. 
This is due to the definition of SF, which is an integral quantity depending on the metabolized drug (see Eq~\eqref{defSF}).
This consideration also extends to the other parameters that influence the surviving fraction, resulting in a greater $\mu^*_i(SF(21600\, s))$ for all $i$ considered in the analysis (Figure \ref{fig:SF_0D}). Lastly, $\sigma_i$ presents minor variations in all parameters, suggesting that interactions between the parameters play a role in the entire time domain.

\subsection{Sensitivity Analysis for 3D-1D-0D model.}

In this section, we systematically adopt the multiscale 3D--1D--0D model introduced in Section \ref{sec:oneway_interaction} for the sensitivity analysis. This choice is motivated by the significant computational cost of the full nonlinear 3D-1D model, which would render comprehensive sensitivity analysis impractical. By replacing the nonlinear metabolic terms with time-dependent surrogate functions obtained from the 0D model, the 3D-1D-0D formulation enables efficient simulation while preserving the key physiological dependencies required to assess parameter influence across scales. {The assumption of steady-state conditions for microvascular flow and oxygen transport (see Sections~\ref{sec:flowmod} and~\ref{OxyMod}) is motivated by the physiological separation of time scales. Vascular adaptation and oxygen equilibration typically occur within seconds to minutes, while prodrug activation and cellular response evolve over hours. This justifies approximating the microvascular and oxygen fields as quasi-steady in the present framework. 
However, neglecting transient fluctuations may obscure the influence of dynamic hypoxia or rapid vascular remodeling on drug activation and cell fate. 
As a consequence, the sensitivity analysis reported here should be interpreted as capturing the dominant parametric trends under averaged conditions, while potentially underestimating the effect of short-lived oxygen fluctuations.}
The sensitivity analysis of the 0D model suggests that the concentration of TPZ in the vasculature is a dominant parameter in determining the concentration of TPZ in the tissue domain.
$\mu^*_{c_{v0}^{tpz}}(c_t^{tpz})$ is approximately two orders of magnitude larger than the secondary influential parameter $\mu^*_{D^{tpz}}(c_t^{tpz})$ (see Figure~\ref{ctaverage}).
In addition to $c_{v0}^{tpz}$, the parameters $k_{met}$, $K$, $\alpha$, $c_{v0}^{ox}$, $V_{max}^{ox}$, and $P_{ox}$ have been identified as relevant for the selected quantities of interest.
For this reason, sensitivity analysis is conducted on this streamlined group of seven parameters for the 3D-1D-0D model. 
The ranges of investigation are prescribed in Table~\ref{tab:boundsSA}.

\begin{figure}[H]
\centering
\subfigure{\textbf{a.}}{\includegraphics[scale=0.25]{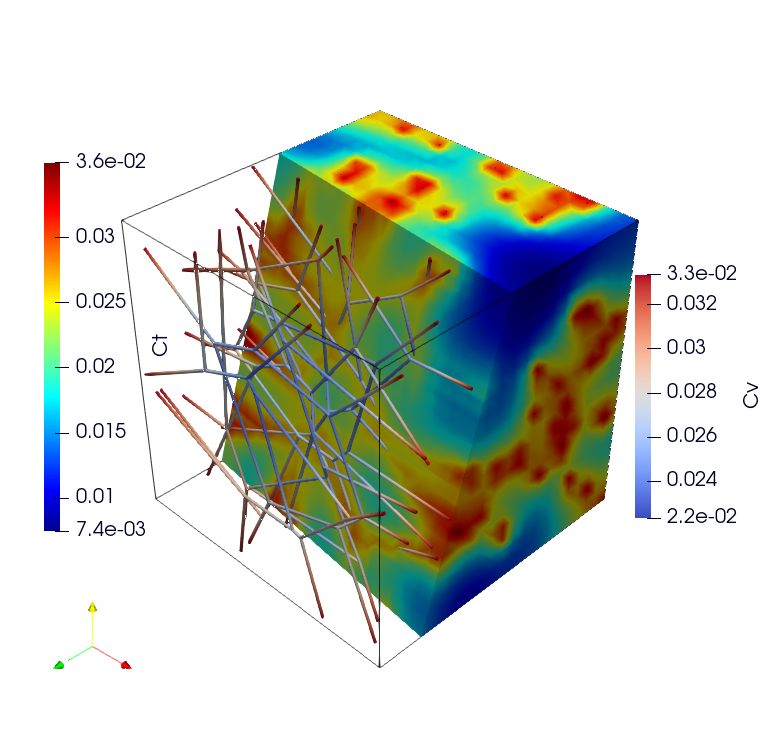}}
\subfigure{\textbf{b.}}{\includegraphics[scale=0.25]{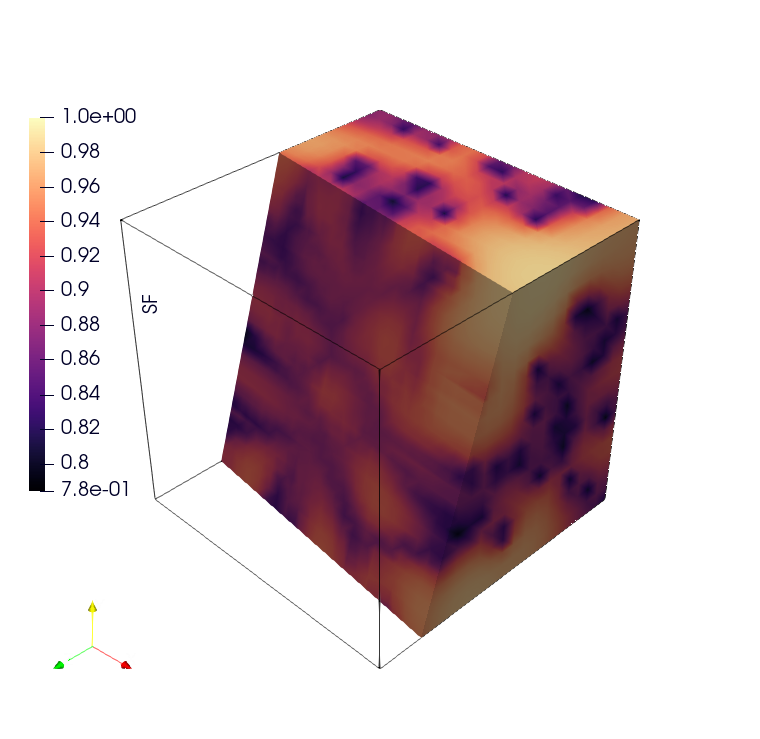}}
\caption{Spatial distribution of $c_t^{tpz}$ along with the distribution of $c_v^{tpz}$ in the vascular network (\textbf{a.}) and $SF$ (\textbf{b.}) taken at $7200\, s$ }
\label{fig:cv_ct_7200}
\end{figure}

As a clear difference among the two approaches, the 3D-1D-0D model provides the spatial distributions of $c_t^{tpz}$ and $SF$, as well as of all other fields involved, as a function of the vascular network immersed (Figure~\ref{fig:cv_ct_7200}.a). 
Specifically, the spatial distribution of TPZ in tissue is a function of vascularization.
Regions with lower drug concentrations present a small number of capillaries immersed (low vascularization), while regions with a large number of capillaries immersed (high vascularization) exhibit a higher drug concentration. 
Analogously, Figure~\ref{fig:cv_ct_7200}.b shows that the spatial distribution of $SF$ is also a function of vascularization. 
To compare the Morris indices associated with the 0D model and those relative to the 3D-1D-0D model, the spatial averages of $c_t^{tpz}$ and $SF$ are calculated and then used as quantities of interest in the sensitivity analysis.

\begin{figure}
\centering
\includegraphics[width=0.8\textwidth]{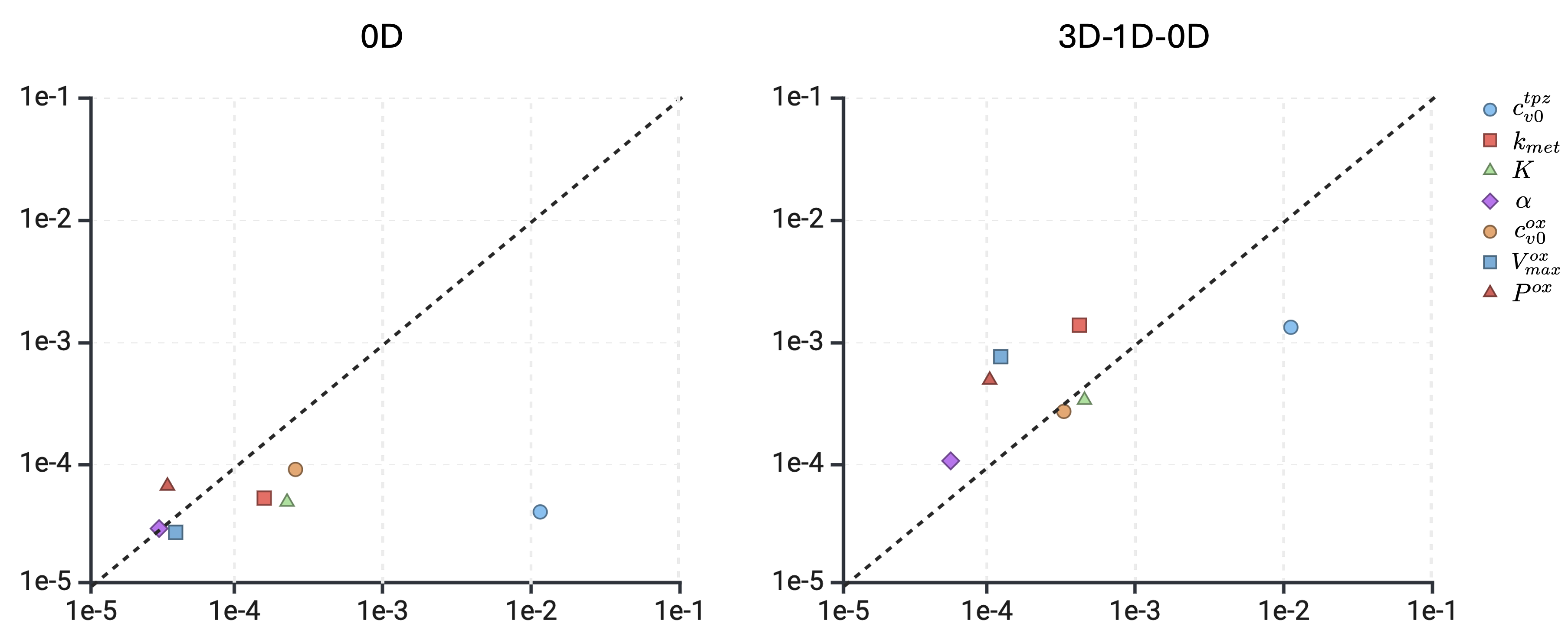}
\caption{Morris indices in the $\mu^*$-$\sigma$ plane (the horizontal axis shows $\mu^*$ and the vertical axis shows $\sigma$) for $c_{v0}^{tpz}$, $k_{met}$, $K$, $\alpha$, $c_{v0}^{ox}$, $V_{max}^{ox}$, and $P_{ox}$ relatively to $\overline{c}_t^{tpz}$.}
\label{fig:3D0D_avct}
\end{figure}

\begin{figure}
    \centering
    \includegraphics[width=0.5\linewidth]{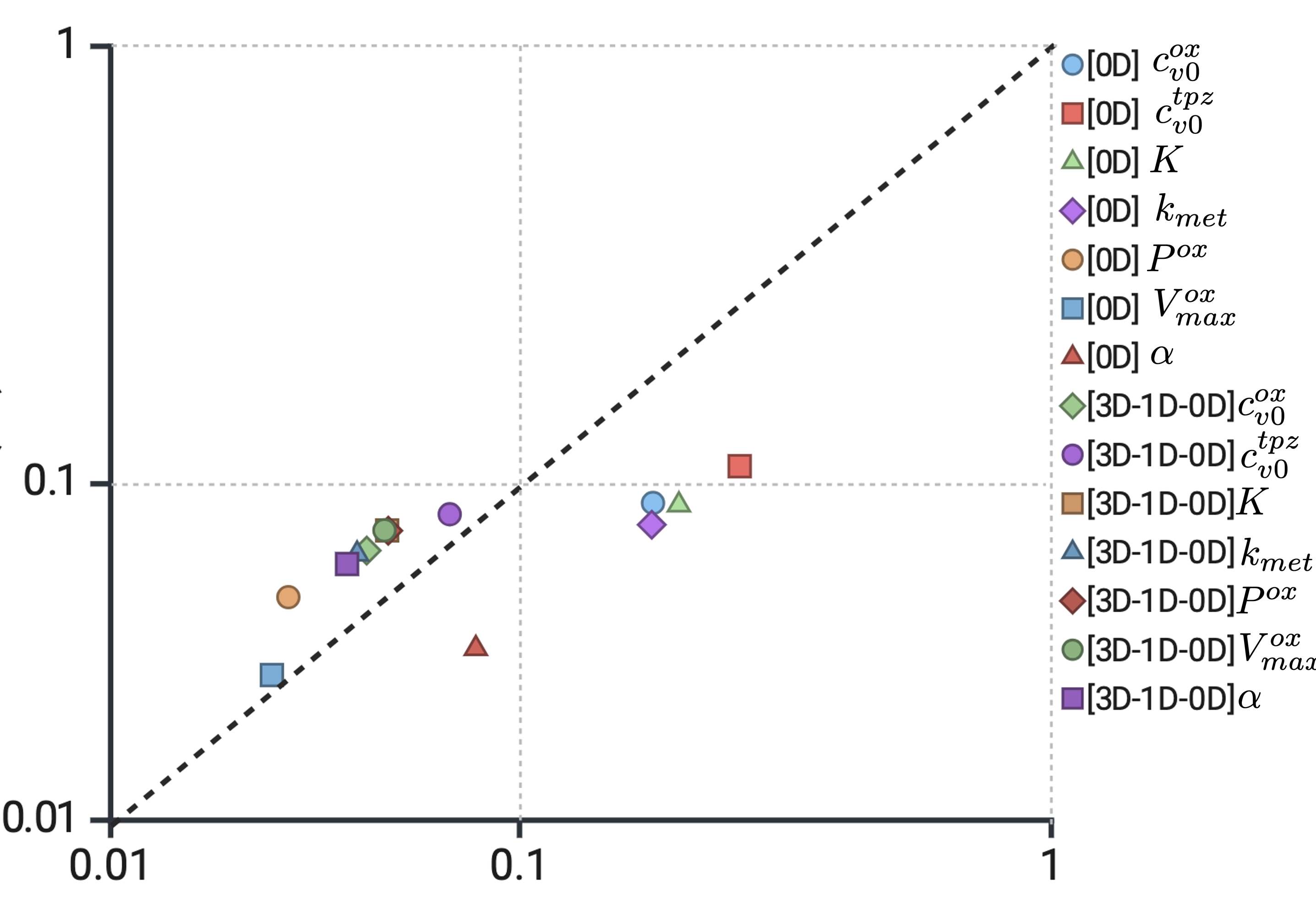}
    \caption{Morris indices in the $\mu^*$-$\sigma$ plane relative to $\overline{SF}$ for the 3D-1D-0D model (the horizontal axis shows $\mu^*$ 
and the vertical axis shows $\sigma$).}
    \label{fig:mu_avsf}
\end{figure}

The Morris indices $\mu^*_{c_{v0}^{tpz}}(\cdot)$ and $\sigma_{c_{v0}^{tpz}}(\cdot)$ with QoI $(\cdot) = \{ \overline{c}_t^{tpz}; \ c_t^{tpz}(7200\, s); \ c_t^{tpz}(10800\, s); \ c_t^{tpz}(21600\, s)\}$ are compared directly to those obtained with 0D and reported in Table~\ref{tab:musigma_ct_3D0D}. 
\begin{table}
\centering
\footnotesize
\begin{tabular}{|c|c|c|c|c|c|}
\hline
 \textbf{index} & \textbf{model} & $\overline{c}_t^{tpz}$ & $c_t^{tpz}(7200\, s)$ & $c_t^{tpz}(10800\, s)$ & $c_t^{tpz}(21600\, s)$\\
\hline
\hline
\multirow{2}{*}{$\mu^*_{c_{v0}^{tpz}}$} & \textbf{0D} & $1.17 \times 10^{-2}$ & $2.94 \times 10^{-2}$ & $2.05 \times 10^{-2}$ & $3.37 \times 10^{-4}$ \\ 
\cline{2-6} 
& \textbf{3D-1D-0D} & $1.13 \times 10^{-2}$ & $2.83 \times 10^{-2}$ &
$1.94 \times 10^{-2}$ & $4.72 \times 10^{-4}$ \\
\hline
\multirow{2}{*}{$\sigma_{c_{v0}^{tpz}}$} & \textbf{0D} & $3.99 \times 10^{-5}$  & $1.73 \times 10^{-4}$ & $6.36 \times 10^{-5}$ &
$4.72 \times 10^{-6}$ \\
\cline{2-6}
& \textbf{3D-1D-0D} & $1.32 \times 10^{-3}$ & $3.30 \times 10^{-3}$ 
& $2.26 \times 10^{-3}$ & $5.51 \times 10^{-5}$  \\ 
\hline
\end{tabular}
\caption{Morris indices for the two model related to $c_{v0}^{tpz}$ for $\overline{c}_t^{tpz}$, $c_t^{tpz}(7200\, s)$, $c_t^{tpz}(10800\, s)$, and $c_t^{tpz}(21600\, s)$.}
\label{tab:musigma_ct_3D0D}
\end{table}
We can observe that the $\mu^*_{c_{v0}^{tpz}}$ indices are very similar between the two models, and the considerations enlightened for the surrogate model are confirmed. 
However, the $\sigma_{c_{v0}^{tpz}}$ indices obtained for the 3D-1D-0D model are approximately two orders of magnitude larger than those of the 0D model, except for $c_t^{tpz}(21600\, s$). 
This is because the 3D-1D-0D model better captures the inherent complexity of the biophysical problem by explicitly including the spatial domain.
Consequently, the Morris indices emphasize the nonlinear relationship between input parameters and output fields. 
As demonstrated in Figure~\ref{fig:3D0D_avct}, this discrepancy is also evident in the $\mu^*-\sigma$ plane for $\overline{c}_t^{tpz}$ relative to the other investigation parameters.
\begin{figure}
\centering
\includegraphics[width=0.33\linewidth]{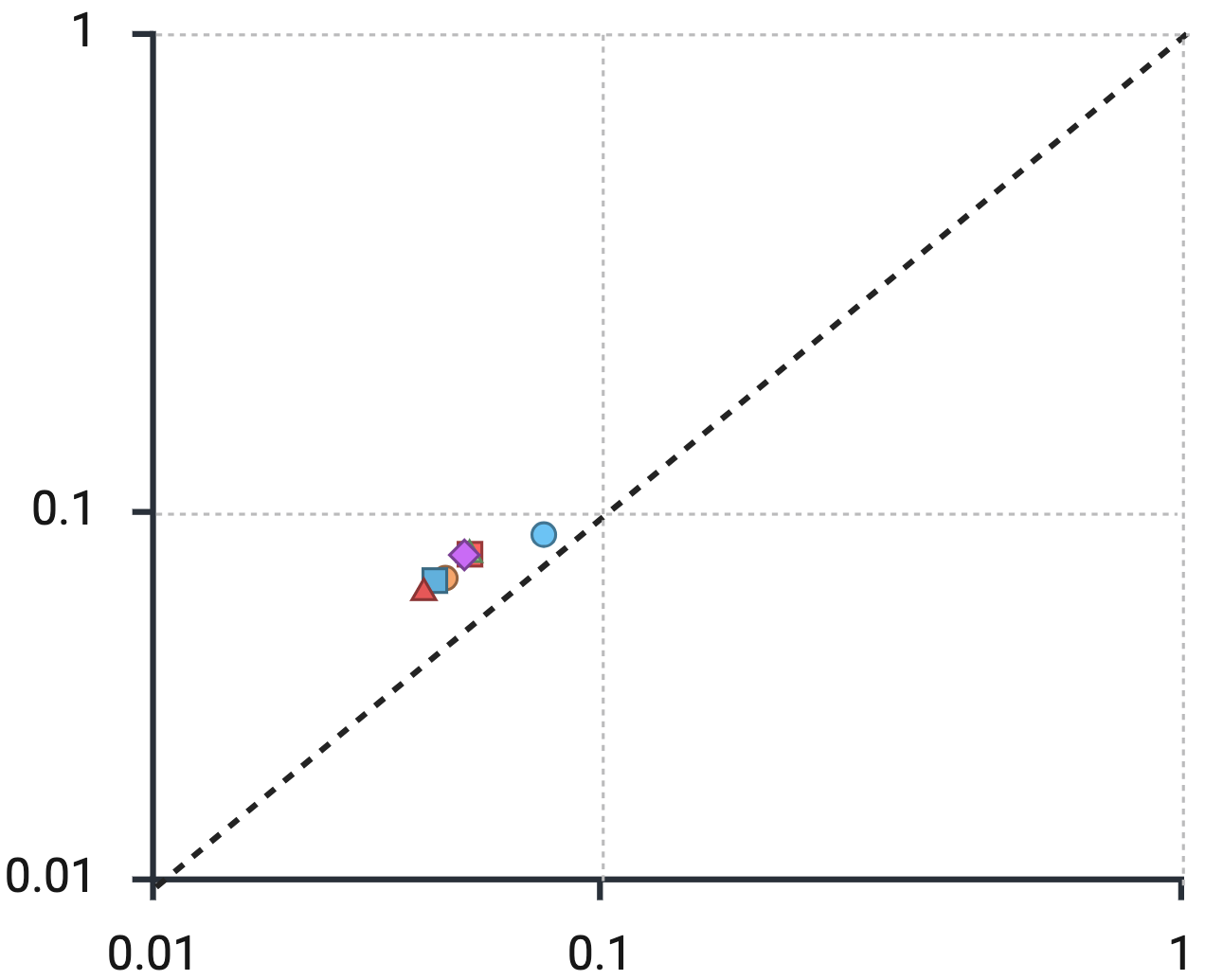}
\hspace{-0.2cm}
\includegraphics[width=0.33\linewidth]{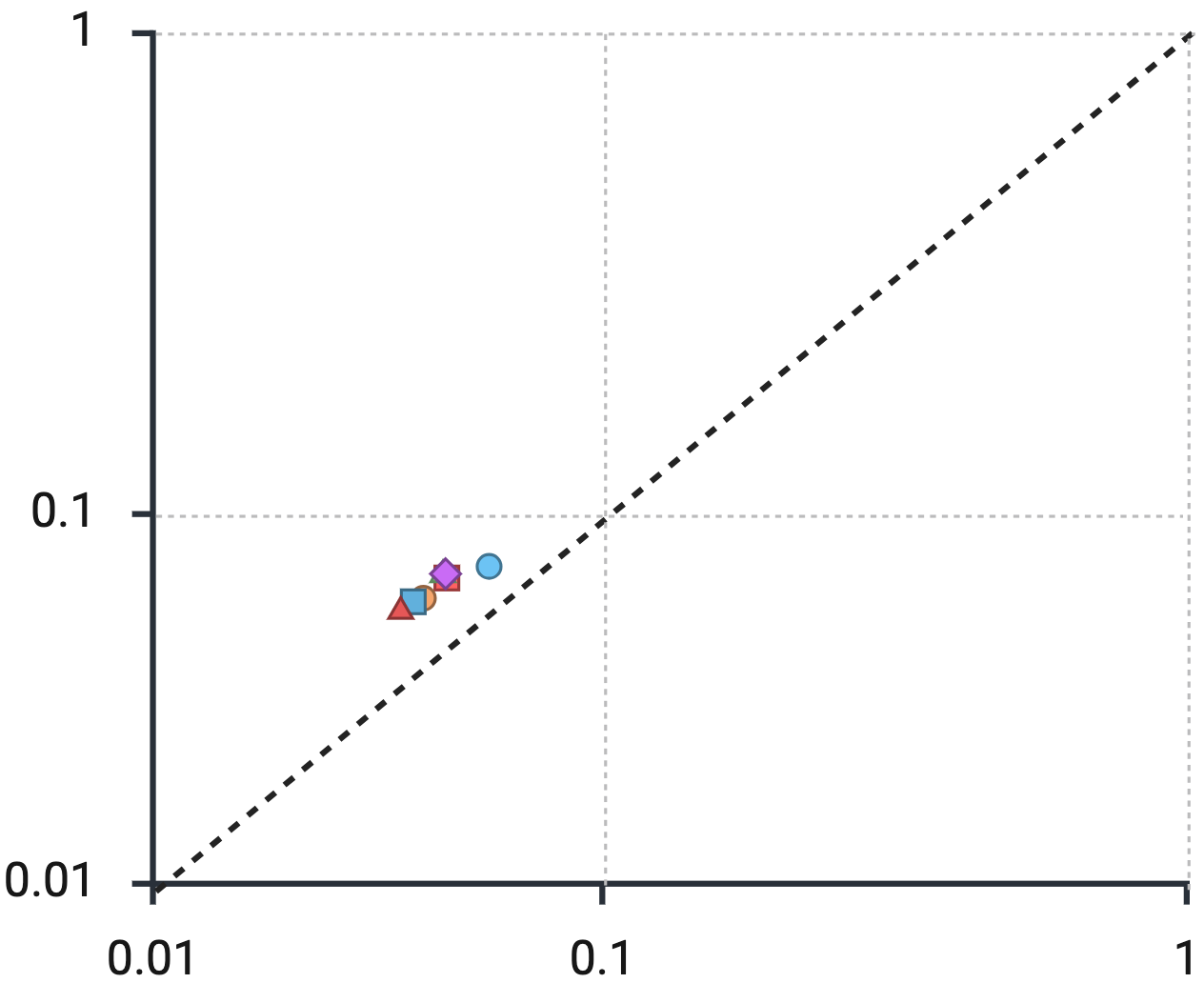}
\hspace{-0.2cm}
\includegraphics[width=0.33\linewidth]{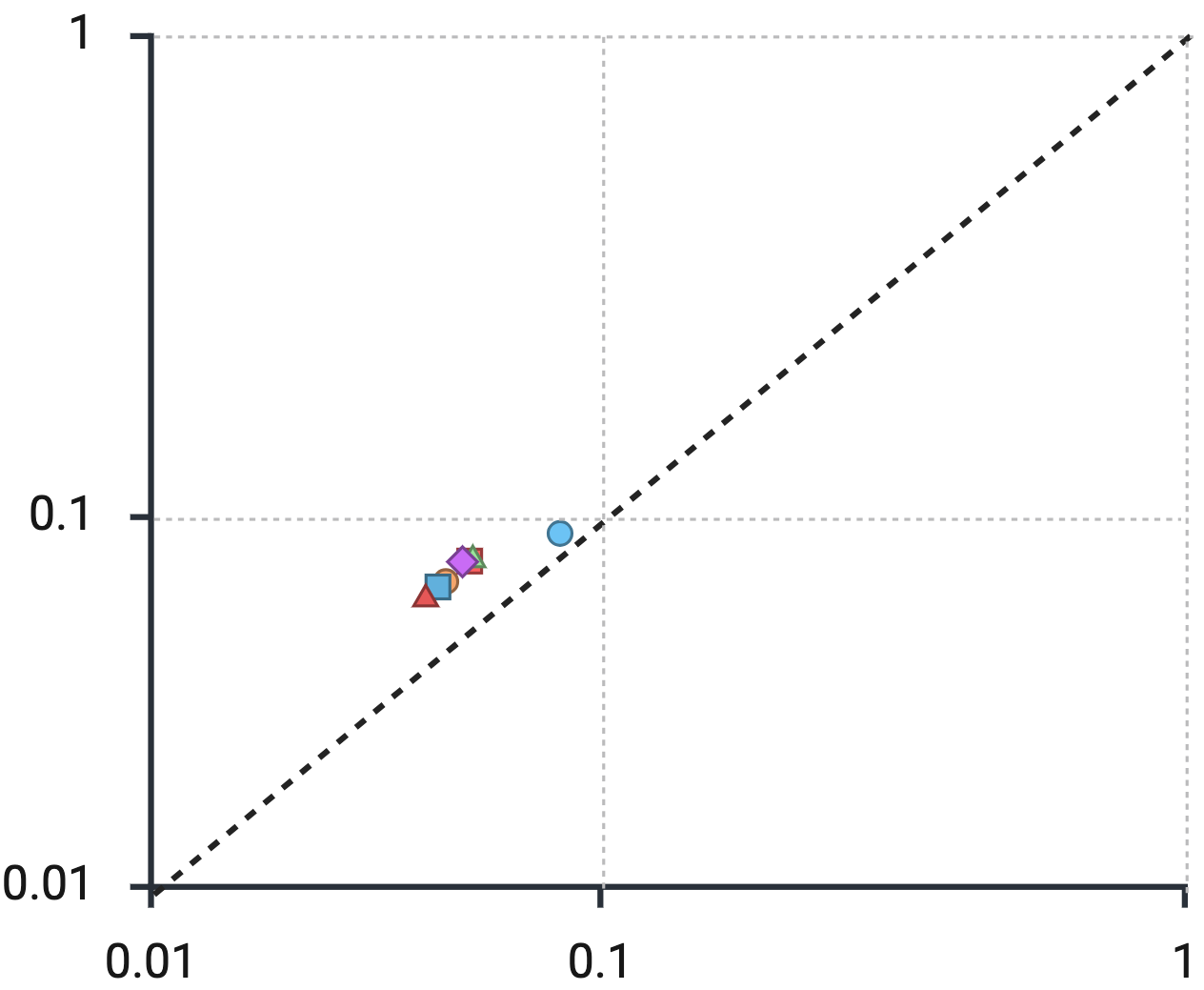} \\
\includegraphics[width=0.5\linewidth]{legenda2_new.png}
\caption{Morris indices in the $\mu^*$-$\sigma$ plane (the horizontal axis shows $\mu^*$ and the vertical axis shows $\sigma$) taken at $7200$, $10800$, and $21600$ seconds, relatively to $SF$, obtained with the 3D-1D-0D-model.}
\label{fig:SF_3D-1D}
\end{figure}
To assess the sensitivity of the surviving fraction computed with the 3D-1D-0D model to the selected input parameters, the Morris indices are first calculated against $\overline{SF}$ (see Figure~\ref{fig:mu_avsf}).  It is noticeable that the indices relative to the 3D-1D-0D model are about one order of magnitude smaller than those of the 0D model, whereas the relative standard deviation is larger. Those obtained larger variances essentially indicate a wider interaction between the parameters and non-linear relations with the averaged distribution of $SF$. In contrast, the smaller $\mu^*_i$ represents a model less sensitive to variations in input parameters. This picture is confirmed by looking at the indices computed for $SF(7200\, s)$, $SF(10800\, s)$, and $SF(21600\, s)$ collected in Figure~\ref{fig:SF_3D-1D}. 

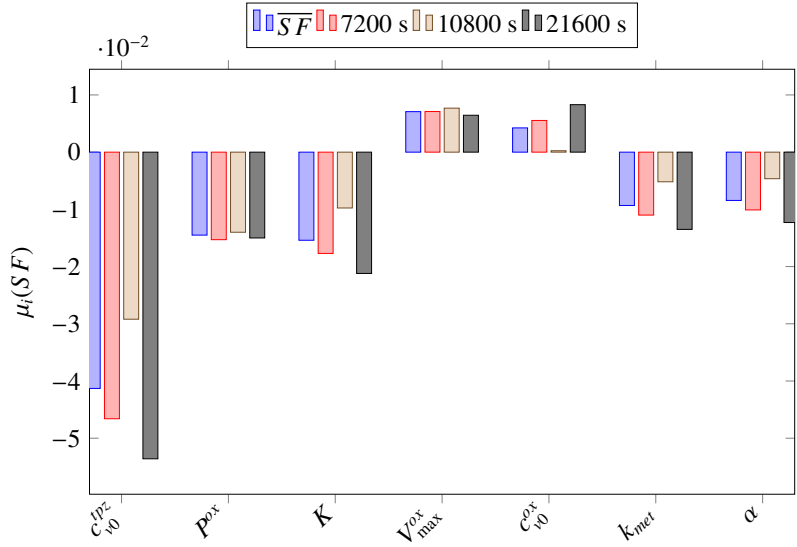
\begin{figure}[H]
\centering
\begin{adjustbox}{scale=0.8}
\begin{tikzpicture}
\begin{axis}[
    ybar,
    bar width=7pt,
    log origin=infty,
    ylabel={$ \mu_i(SF)$},
    symbolic x coords={
        $c_{v0}^{tpz}$, $P^{ox}$, $K$, $V_{\text{max}}^{ox}$,
        $c_{v0}^{ox}$, $k_{met}$, $\alpha$
    },
    xtick=data,
    x tick label style={rotate=45, anchor=east},
    enlarge x limits=0.05,
    x=50,
    bar shift auto,
    legend style={at={(0.5,1.05)}, anchor=south, legend columns=-1},]

% Mean SF
\addplot+ coordinates {
    ($c_{v0}^{tpz}$,-0.0413) ($P^{ox}$,-0.0145) ($K$,-0.0154)
    ($V_{\text{max}}^{ox}$,0.00708) ($c_{v0}^{ox}$,0.00424)
    ($k_{met}$,-0.00932) ($\alpha$,-0.00844)
};
\addlegendentry{$\overline{SF}$}

% SF(7200 s)
\addplot+ coordinates {
    ($c_{v0}^{tpz}$,-0.0466) ($P^{ox}$,-0.0153) ($K$,-0.0177)
    ($V_{\text{max}}^{ox}$,0.00710) ($c_{v0}^{ox}$,0.00554)
    ($k_{met}$,-0.0110) ($\alpha$,-0.0101)
};
\addlegendentry{7200 s}

% SF(10800 s)
\addplot+ coordinates {
    ($c_{v0}^{tpz}$,-0.0292) ($P^{ox}$,-0.0140) ($K$,-0.00976)
    ($V_{\text{max}}^{ox}$,0.00769) ($c_{v0}^{ox}$,0.000249)
    ($k_{met}$,-0.00518) ($\alpha$,-0.00464)
};
\addlegendentry{10800 s}

% SF(21600 s)
\addplot+ coordinates {
    ($c_{v0}^{tpz}$,-0.0536) ($P^{ox}$,-0.0150) ($K$,-0.0212)
    ($V_{\text{max}}^{ox}$,0.00645) ($c_{v0}^{ox}$,0.00830)
    ($k_{met}$,-0.0135) ($\alpha$,-0.0123)
};
\addlegendentry{21600 s}

\end{axis}
\end{tikzpicture}
\end{adjustbox}
\caption{Bar chart of $\mu_i(\overline{SF})$ and $\mu_i(SF)$ taken at $7200$, $10800$, and $21600$ obtained with the 3D-1D-0D-model.}
\label{fig:musign_3D-1D}
\end{figure}

Until now, our analysis has been based primarily on the indices \(\mu^*_i\), defined as the mean of the absolute values of the elementary effects (EE). This choice provides a robust estimate of the overall importance of each input parameter in the output, regardless of the direction of influence. However, additional insight can be gained by also examining the signed mean \(\mu_i\) of the elementary effects, which preserves the directionality of the input-output relationship. While \(\mu^*_i\) is suitable for identifying which parameters have strong effects, the sign of \(\mu_i\) helps to interpret whether an increase in a given parameter tends to increase or decrease the output, on average.
{In the Morris screening method, the mean of the absolute elementary effects, $\mu_i^*$, is the standard metric to quantify the overall importance of a factor, as it avoids the cancellation of opposing contributions. For this reason, in our analysis, we consider $\mu_i^*$ as the primary indicator of parameter relevance. Alongside $\mu_i^*$, we also report the signed mean $\mu_i$, not as a measure of global importance but rather to provide qualitative information on the direction of the effect. This distinction allows us to combine the robustness of $\mu_i^*$ in ranking influential factors with the interpretability of $\mu_i$ in indicating the prevailing sign of the effect.}
We apply this directional analysis to the surviving fraction \(SF\), a key indicator of therapeutic efficacy. 
Figure~\ref{fig:musign_3D-1D} displays the signed indices \(\mu_i(\overline{SF})\), \(\mu_i(SF(7200\,\mathrm{s}))\), \(\mu_i(SF(10800\,\mathrm{s}))\), and \(\mu_i(SF(21600\,\mathrm{s}))\) for the seven parameters previously identified as the most influential. 
Negative values of \(\mu_i\) suggest that increasing the corresponding parameter tends to reduce the surviving fraction, i.e., enhance the cytotoxic effect of TPZ. 
In contrast, positive values indicate that increases in the parameter are associated with a reduction in drug efficacy.
Hence, they should be interpreted alongside \(\mu^*_i\) and \(\sigma_i\) to provide a more complete understanding of sensitivity in the model.
Specifically, we observe negative values of \(\mu_i\) for \(i = \{c_{v0}^{tpz}, P^{ox}, K, k_{met}, \alpha\}\), indicating that increases in these parameters tend to decrease \(\overline{SF}\). 
The effect of \(c_{v0}^{tpz}, K, k_{met}, \alpha\) aligns with biological intuition: all these parameters are directly involved in drug availability or metabolism and thus influence TPZ activation. 
The role of \(P^{ox}\), \(c_{v0}^{ox}\), and \(V_{max}^{ox}\) is more nuanced. 
Surely, an increased vascular oxygen permeability \(P^{ox}\) improves tissue oxygenation, reducing TPZ consumption \(m^{tpz}\).
In our data, an increased \(P^{ox}\) generally results in a lower \(SF\).
Conversely, even if a greater \(c_{v0}^{ox}\) does increase the tissue oxygenation, the general result is a higher \(SF\).
A similar effect on \(SF\) is obtained when rising \(V_{max}^{ox}\), which is expected to decrease the tissue oxygenation.
In summary, these parameters affect tissue oxygenation and modify TPZ consumption.
Consequently, \(c_t^{tpz}\) increases, leading to TPZ accumulation within the tissue.
Following the definition by the equation \ref{defSF}, the \(SF\) is affected by both the TPZ concentration ($c_t^{tpz}$) and the drug metabolization \(m^{tpz}\).
Varying these parameters, we obtain a decrease in drug consumption but an increase in drug concentration, resulting in a behavior hardly predictable a priori.
Our modelling approach enables the estimation of \(SF\) as a combination of these two effects, accounting for non-linear interactions due to the inherent coupling of these phenomena within the microenvironment.
We remark that the definition \(SF\) (Eq.~\ref{defSF}) could also be refined and substituted in the model based on further experimental evidence, separating the contributions related to drug metabolization and drug concentration.

\section{Conclusions and Future Developments}

This study presents a mathematical and mechanistic framework for exploring the pharmacokinetics and pharmacodynamics of hypoxia-activated drugs in solid tumors. By integrating a spatially resolved 3D-1D model of vascular flow and oxygen transport with a 0D surrogate model that describes nonlinear drug metabolism and cell survival, we provide a computational approach to investigate how physical, physiological, and geometric factors in the tumor microenvironment influence drug efficacy.
This multiscale model captures key feedback mechanisms between drug activation and local oxygen concentration, reflecting the key challenges associated with the delivery of hypoxia-activated drugs, such as Tirapazamine. Through global sensitivity analysis, we systematically identified the parameters that most influence drug distribution and therapeutic outcome, highlighting the dominant role of vascular drug concentration, oxygen availability, and metabolic rates.

Beyond methodological contributions, our findings have important implications for the clinical optimization of TPZ and other hypoxia-activated prodrugs. The multiscale 3D-1D-0D framework developed in this work provides a computationally efficient yet physiologically grounded tool to simulate the spatiotemporal evolution of drug concentration and therapeutic response within the hypoxic tumor microenvironment. By integrating nonlinear pharmacodynamic effects via exogenous surrogate functions, the model enables rapid sensitivity analyses and parametric studies that would otherwise be prohibitive with full-scale simulations.
This capability opens the door to systematic exploration of treatment protocols, such as timing, dosage, and vascular delivery of TPZ, under varying levels of tissue oxygenation. Furthermore, the surrogate-based architecture supports future incorporation of patient-specific vascular geometries and clinical imaging data, paving the way for predictive simulation platforms that help customize hypoxia-targeted therapies. As such, the proposed modeling approach may serve as a foundation for digital twin frameworks aimed at optimizing the therapeutic index of bioreductive agents in precision oncology.

Although the proposed model offers a detailed representation of the tumor microenvironment on multiple scales, the results should primarily be interpreted as exploratory. The analysis demonstrates internal consistency and physiological plausibility in a range of hypothetical scenarios, but it does not yet include calibration to specific experimental data sets or patient-derived measurements. Therefore, its primary value lies in the generation of hypotheses and the guidance of future experimental or computational studies, rather than providing predictive assessments for clinical settings.
Furthermore, it has several limitations that suggest opportunities for future improvements. Some biological processes were simplified or excluded, including dynamic changes in vessel permeability, red blood cell interactions, vascular remodeling, and the time-dependent effects of therapy on the vasculature. {This assumption implies that acute vascular remodeling and transient fluctuations in oxygen supply are neglected. As a result, the predicted oxygen distributions are smoother and more stable than \textit{in vivo}, and the corresponding survival fraction $SF$ reflects the average conditions rather than the effects of cycling hypoxia or abrupt vascular alterations. Although this simplification is sufficient for our sensitivity analysis, it may underestimate the impact of transient hypoxic events on drug activation and cellular response.} Furthermore, steady-state or quasi-steady oxygen transport is assumed, and the model lacks a complete temporal coupling between all components, which could be important for modeling rapid response to treatment or combination therapies. {In addition to highlighting the exploratory character of this work, the sensitivity analysis suggests several testable hypotheses that can guide future experimental investigations. 
First, increasing vessel oxygen permeability enhances tissue oxygenation but simultaneously decreases hypoxia-driven drug activation, leading to a non-monotonic impact on cell survival. 
Second, variability in metabolic parameters $(V_{\max}^{\mathrm{ox}}, K_m^{\mathrm{ox}})$ across cell lines or in patient-specific settings is predicted to strongly influence therapeutic outcome, suggesting that these quantities are potential biomarkers of treatment response. Third, our results indicate that the sensitivity of the survival fraction to microvascular geometry is secondary compared to systemic pharmacokinetics and metabolic rates, implying that optimizing drug scheduling may outweigh the effects of structural vascular heterogeneity. By making these hypotheses explicit, the model can serve not only as a computational tool but also as a generator of experimentally verifiable predictions.}

Future work will address these aspects by incorporating time-dependent vascular remodeling, angiogenesis, and tumor growth models, as well as refining drug-specific metabolic pathways. Further sensitivity analysis of geometrical parameters and their interactions could yield more insight into patient-specific variability. Integration with \textit{in vivo} and \textit{in vitro} experimental data will be crucial to calibrate the model and support uncertainty quantification, inverse modeling, and personalized therapy design. Furthermore, integrating the current framework with data-driven techniques could offer practical strategies for real-time adaptation to treatment. {On the other hand, future work will explore extensions of the present model, including fractional-order formulations, to capture memory effects and anomalous transport phenomena that may further enhance predictive capability.}

In summary, this work provides a mechanistic tool for probing the complex interactions between vascular architecture, oxygen dynamics, and drug metabolism in the context of hypoxia-targeted therapies. By identifying influential parameters and mechanistic drivers of treatment response, it lays the foundations for more comprehensive studies aimed at optimizing therapeutic strategies and improving our understanding of treatment resistance in hypoxic tumors. {Moreover, our results highlight the crucial role of oxygen heterogeneity in shaping the activation and efficacy of hypoxia-activated prodrugs, demonstrating the advantages of a multiscale 3D-1D-0D approach over purely surrogate or fully spatial models. Looking ahead, the framework can be extended to incorporate time-dependent vascular dynamics, integrated with experimental datasets for quantitative validation, and enriched with additional biological processes such as DNA damage repair or drug resistance. 
These developments will further enhance the predictive power of the model and support its application for treatment optimization in clinically relevant settings. Although the present framework is not yet intended for predictive assessments in clinical settings, it naturally lends itself to progressive refinement toward clinical applicability. A possible roadmap involves: (i) preclinical validation against \textit{in vitro} and \textit{in vivo} data, calibrating oxygen transport and pharmacokinetics/pharmacodynamics parameters; (ii) incorporation of patient-specific information, such as imaging-derived hypoxia maps and pharmacokinetics profiles, to personalize inputs; (iii) systematic integration of biological variability, to provide predictive ranges rather than single-point estimates; and finally  (iv) deployment in \textit{in silico} clinical trials to explore dosing and scheduling strategies before translation into actual patient studies.}

\section*{Author contributions}
\textbf{A. Coclite, R. Montanelli Eccher, L. Possenti, P. Vitullo}: Conceptualization, Methodology,  Writing — review and editing, Writing — original draft preparation, Supervision.
\textbf{P. Zunino}: Conceptualization, Methodology, Funding and resources acquisition, Writing — review and editing, Writing — original draft preparation, Supervision.

\section*{Acknowledgments}
 AC acknowledges support from the MUR PRIN 2022 project (code 2022M9BKBC, CUP D53D23005880006) and the MUR PRIN PNRR 2023 project (code P202254HT8, CUP B53D23027760001). PV and PZ acknowledge support from the MUR PRIN 2022 project 2022WKWZA8 \textit{Immersed methods for multiscale and multiphysics problems} (IMMEDIATE), part of the Next Generation EU program Mission 4, Component 2, CUP D53D23006010006. This work has received funding from the European Union’s Euratom research and training programme 2021–2027 under grant agreement No. 101166699 \textit{Radiation safety through biological extended models and digital twins} (TETRIS). This work was supported by the ERA-PerMed initiative through the Joint Transnational Call 2024, under the project \textit{High-resolution, risk-based, and outcome-driven digital twins for precision radiotherapy in breast cancer} (Hi-ROC), project ID: EP PerMed-JTC2024-187, co-funded by the European Union under the Horizon Europe Framework Programme and national funding agencies. The work was also supported by the AIRC Investigator Grant no. IG21479 (PI: Tiziana Rancati). This research is part of the activities of the \textit{Dipartimento di Eccellenza} 2023--2027, Department of Mathematics, Politecnico di Milano.  AC, PV, and PZ are members of the Gruppo Nazionale per il Calcolo Scientifico (GNCS) of the Istituto Nazionale di Alta Matematica (INdAM).
 
\section*{Conflict of interest}
The authors declare that they have no known competing financial interests or personal relationships that could have appeared to influence the work reported in this paper.

\appendix

\section{Nomenclature}
\label{app-A}

\begin{itemize}
  \item \textbf{TPZ} - Tirapazamine, a representative hypoxia-activated prodrug.
  \item \textbf{SF} - Surviving Fraction, quantifying the proportion of viable cells over time.
  \item \textbf{TME} - Tumor Microenvironment.
  \item \textbf{RBC} - Red Blood Cell.
  \item \textbf{3D--1D model} - Multiscale model coupling three-dimensional tissue transport with one-dimensional vascular networks.
  \item \textbf{0D model} - Surrogate pharmacokinetics–pharmacodynamics model used to simulate average tissue-level drug and oxygen dynamics.
  \item \textbf{3D--1D--0D model} - Multiscale model combining spatially resolved transport (3D–1D) with surrogate terms derived from the 0D model.
  \item \textbf{PK/PD} - Pharmacokinetics/Pharmacodynamics.
  \item \textbf{SA} - Sensitivity Analysis.
  \item \textbf{DoE} - Design of Experiments.
  \item \textbf{EE} - Elementary Effects method, used for global sensitivity analysis.
  \item \textbf{KS test} - Kolmogorov–Smirnov test, used to assess distributional similarity.
  \item \textbf{ODE} - Ordinary Differential Equation.
  \item \textbf{PDE} - Partial Differential Equation.
  \item \textbf{FEM} - Finite Element Method.
  \item \textbf{DOF} - Degrees of Freedom.
\end{itemize}

\section{List of Symbols}\label{app-B}

\begin{itemize}
  \item \( \Omega \) - 3D tissue domain [\(\mathrm{mm}^3\)]
  \item \( \Lambda \) - 1D vascular network domain [\(\mathrm{mm}\)]
  \item \( \Gamma \) - Vessel-tissue interface [\(\mathrm{mm}^2\)]
  \item \( \mathbf{x} \in \Omega \) - Spatial position in the tissue domain
  \item \( s \in \Lambda \) - Arc-length parameter along 1D vessel centerlines
  \item \( t \) — Time [s]
  \item \( c_t^{tpz}(\mathbf{x}, t) \) - TPZ concentration in tissue [mol/m\(^3\)]
  \item \( c_v^{\mathrm{tpz}}(s, t) \) - TPZ concentration in vessels [mol/m\(^3\)]
  \item \( c_t^{ox}(\mathbf{x}, t) \) - Oxygen concentration in tissue [mol/m\(^3\)]
  \item \( c_v^{\mathrm{ox}}(s, t) \) - Oxygen concentration in vessels [mol/m\(^3\)]
  \item \( \phi \) - Population of viable cells
  \item \( P \) - Permeability coefficient across the vessel wall [m/s]
  \item \( S/V \) - Surface-to-volume ratio of the vessel–tissue interface [1/m]
  \item \( SF(t) \) - Surviving fraction at time \( t \) (dimensionless)
  \item \( \mathbf{u} \) - Interstitial fluid velocity field [m/s]
  \item \( D_t^{ox}, D_v^{ox} \) - Diffusion coefficients for oxygen [m\(^2\)/s]
  \item \( D_t^{tpz}, D_v^{tpz} \) - Diffusion coefficients for TPZ [m\(^2\)/s]
  \item \( \mathcal{M}_{ox}(c^{ox}, SF) \) - Oxygen consumption rate [mol/(m\(^3\)·s)]
  \item \( \mathcal{M}_{tpz}(c_t^{tpz}, c_t^{ox}, SF) \) - TPZ metabolism rate [mol/(m\(^3\)·s)]
  \item \( m^{tpz} \) - Surrogate TPZ metabolic rate [s\(^{-1}\)]
  \item \( m^{ox} \) - Surrogate oxygen consumption rate [s\(^{-1}\)]

  \item \( r(t) \) - Effective metabolic efficiency at time $t$ [s\(^{-1}\)]
  \item \( K \) - Half-saturation constant for TPZ oxygen modulation [mol/m\(^3\)]
  \item \( k_{met} \) - Linear metabolism rate constant [s\(^{-1}\)]
  \item \( V_{max}^{tpz} \) - Maximum velocity of nonlinear metabolism [mol/(m\(^3\)·s)]
  \item \( K_m \) - Michaelis–Menten constant for TPZ [mol/m\(^3\)]
\end{itemize}

\end{document}